\documentclass[aos,preprint]{imsart}
\usepackage{enumerate}
\usepackage{amsmath}
\usepackage{amssymb}
\usepackage{amsthm}
\renewcommand{\qed}{\alignleft {\text{\tiny{$\blacksquare$}}}}
\usepackage{amsfonts}
\usepackage{graphicx}

\theoremstyle{plain}
\newtheorem{theorem}{Theorem}[section]	
\newtheorem{lemma}{Lemma}[section]
\newtheorem{corollary}{Corollary}[section]

\theoremstyle{definition}
\newtheorem{definition}{Definition}[section]
\newtheorem{remark}{Remark}[section]
\newtheorem{example}{Example}[section]

\newcommand{\R}{\mathbb{R}}

\renewcommand{\qed}{\hfill{\tiny \ensuremath{\blacksquare} }}%

\renewcommand{\Pr}{{\mathrm{P}}}

\usepackage{geometry}

\usepackage{geometry}                
\RequirePackage[OT1]{fontenc}
\RequirePackage[colorlinks]{hyperref}
\newcommand{\DD}{{\rm D}}
\newcommand{\QQ}{{\rm Q}}
\newcommand{\RR}{{\rm R}}
\newcommand{\BL}{{\rm BL}}
\newcommand{\KR}{{\rm BL}}
\newcommand{\E}{{\rm I}\kern-0.18em{\rm E}}
\renewcommand{\P}{{\rm I}\kern-0.18em{\rm P}}
\newcommand{\D}{{\rm I}\kern-0.18em{\rm D}}
\renewcommand{\Pr}{{\rm I}\kern-0.18em{\rm P}}
\renewcommand{\R}{{\rm I}\kern-0.18em{\rm R}}
\newcommand{\B}{{\rm I}\kern-0.18em{\rm B}}
\renewcommand{\phi}{\varphi}
\newcommand{\eps}{\varepsilon}
\newcommand{\mU}{\mathcal{U}}
\newcommand{\mY}{\mathcal{Y}}
\newcommand{\mD}{\mathcal{D}}
\newcommand{\mX}{\mathcal{X}}
\newcommand{\mK}{\mathcal{K}}
\newcommand{\bK}{\mathbb{K}}

\newtheorem*{lemma*}{Lemma}
\newtheorem*{proposition*}{Proposition}

\newtheorem*{definition*}{Definition}
\begin{document}
\begin{frontmatter}

 \title{Monge-Kantorovich depth, \\quantiles, ranks, and signs}
\runtitle{Monge-Kantorovich Depth}


\begin{aug}
\author{\fnms{Victor} \snm{Chernozhukov}\thanks{Research  supported by the NSF Grant SES 1061841.}, \ead[label=e1]{vchern@mit.edu}}
\author{\fnms{Alfred} \snm{Galichon}\thanks{Research   funded by the European
Research Council under the European Union’s Seventh Framework Programme (FP7/2007-2013), ERC grant
agreement \#313699.}, \ead[label=e2]{ag133@nyu.edu}}\\
\author{\fnms{Marc} \snm{Hallin}\thanks{
Research  supported by the IAP research network grant~P7/06 of the Belgian government (Belgian Science
Policy), a Cr\' edit aux Chercheurs of the Fonds National de la Recherche Scientifique,  and the Discovery grant DP150100210 of the Australian Research Council.}, \ead[label=e3]{mhallin@ulb.ac.be}}
\and
\author{\fnms{Marc} \snm{Henry}\thanks{Research  supported by the  SSHRC Grant 435-2013-0292 and the NSERC Grant
356491-2013.} \ead[label=e4]{marc.henry@psu.edu}}

\runauthor{Chernozhukov, Galichon, Hallin, Henry}

\affiliation{MIT, NYU, Universit\' e libre de Bruxelles and Princeton University, Penn State}

\address{\footnotesize
Department of Economics and Center for Statistics \\
Massachusetts Institute of Technology \\
Cambridge, MA 02139, USA\\
\printead{e1}
 }

\address{\footnotesize 
Economics Department and Courant Institute,\\
New York University\\
New York, NY 10012, USA\\  and\\
Sciences Po, Economics Department \\
75007 Paris, France\\ 
\printead{e2}}

\address{ \footnotesize ECARES \\
Universit\' e libre de Bruxelles\\
CP 114 Brussels,
Belgium\\ and\\
ORFE\\
Princeton University\\
Princeton, NJ 08540, USA\\ 
\printead{e3}  }

\address{ \footnotesize Department of Economics\\
The Pennsylvania State University\\
University Park, PA 16802, USA\\
\printead{e4}
}

\end{aug}

\begin{abstract}

We propose new concepts of statistical depth, multivariate quantiles, vector quantiles and ranks, ranks, and signs, based on canonical transportation maps between a distribution of interest on $\R^d$ and a reference distribution on the $d$-dimensional unit ball. The new depth concept, called {\em Monge-Kantorovich depth}, specializes to halfspace depth for $d=1$ and in the case of spherical distributions, but,  for more general distributions, differs from  the latter  in the ability for its contours  to account for non convex features of the distribution of interest. We propose empirical counterparts to the population versions of those Monge-Kantorovich depth contours, quantiles, ranks, signs, and vector quantiles and ranks, and show their consistency  by establishing a uniform convergence property for empirical (forward and reverse) transport maps, which is the main theoretical result  of this paper.
 \end{abstract}



\end{frontmatter}

\section{Introduction} \label{intro}
The concept of statistical depth was introduced in order to overcome the lack of a canonical ordering in $\R^d$ for $d>1$, hence the absence of the related notions of quantile and distribution functions, ranks, and signs.
The earliest and most popular depth concept is  halfspace depth, the definition of which goes back  to Tukey \cite{Tukey:75}. Since then, many other concepts have been considered:~{simplicial depth} \cite{Liu:90},  majority depth (\cite{Singh:91} and \cite{LS:93}), {projection depth} (\cite{Liu:92}, building on \cite{Stahel:81} and \cite{Donoho:82}, \cite{Zuo:2003}), {Mahalanobis depth} (\cite{Mahalanobis:36}, \cite{Liu:92}, \cite{LS:93}),  {Oja depth} \cite{Oja:83},  {zonoid depth} (\cite{KM:97} and  \cite{Koshevoy:2002}), {spatial depth} (\cite{DK:92}, \cite{MO:95}, \cite{Chaudhuri:96}, \cite{VZ:2000}), $L^p$~depth \cite{ZS:2000}, among many others.
An axiomatic approach, aiming at unifying all those concepts, was initiated by Liu \cite{Liu:90} and Zuo and Serfling \cite{ZS:2000}, who list four properties that are generally considered desirable for any statistical depth function, namely affine invariance, maximality at the center, linear monotonicity relative to the deepest points, and vanishing at infinity (see Section~\ref{sec:LZS} for details). Halfspace depth 
  is the prototype of a depth concept satisfying the Liu-Zuo-Serfling axioms for the family~$\mathcal P$ of all absolutely continuous distributions on $\R^d$.

An important feature of halfspace depth is the convexity of its contours, which thus satisfy the star-convexity requirement embodied in the linear monotonicity axiom. That feature  is shared by most existing depth concepts and might be considered undesirable for distributions with non convex supports or level contours, and multimodal ones. Proposals have been made, under the name  of { local depths}, to deal with this, while retaining the spirit of the Liu-Zuo-Serfling axioms: see \cite{CDPB:2009}, \cite{HKV:2010},  \cite{AR:2011}, and \cite{PB:2013} who provide an in-depth discussion of those various attempts. In this paper, we take  a totally different and more agnostic  approach, on the model of the discussion by Serfling in \cite{Serfling:2002}: if the ultimate purpose of statistical depth is to provide, for each distribution $ P$,  a  $ P$-related   ordering   of~$\R^d$  producing adequate  concepts of quantile and distribution functions, ranks and signs, the relevance of  a given depth function should be evaluated in terms of the relevance of the resulting ordering, and the  quantiles, ranks and signs it produces.

Now, the concepts of quantiles, ranks and signs are well understood in two particular cases, essentially,  that should serve as benchmarks. The first case is that of the family~$\mathcal{P}^1$  of all distributions with nonvanishing Lebesgue densities over convex support sets. Here,  the concepts of quantile and distribution functions, ranks, and signs are related to the ``classical" univariate ones.  The second case is that of the family~$\mathcal{P}^d_{\mbox{\scriptsize ell}}$
   of all full-rank  elliptical distributions over~$\R^d$ ($d>1$) with radial densities over elliptical support sets. Recall that the family~$\mathcal{P}^d_{\mbox{\scriptsize ell}; g}=\{P_{\mu ,\Sigma , g}\}$ of elliptical distributions with given { radial density~$g$ and distribution function $G$}  is a parametric family indexed by a location parameter~$\mu$ and a { scatter}  parameter $\Sigma$ (a symmetric positive definite real matrix) such that a random vector $X$ has  distribution $P_{\mu ,\Sigma , g}$ iff the residual $Y:= \Sigma^{-1/2}(X-\mu)$, which results from transforming $X$ into isotropic position, has spherical distribution $P_{0 ,I , g}$. Further, this is equivalent to $\RR_P (Y) = (Y/\|Y\|)G(\|Y\|)$ having the spherical uniform distribution $U_d$ on the unit ball $\mathbb{S}^d$ in $\R^d$.  By {\it spherical uniform}, we mean the distribution of a random vector $r\phi$, where $r$ is uniform on $[0,1]$, $\phi$ is uniform on the unit sphere $\mathcal S^{d-1}$, and~$r$ and $\phi$ are mutually independent.   There, spherical contours with~$ P_{\mu ,I , g}$-probability contents $\tau$ coincide with the  halfspace depth contours, and  provide a natural definition of { $\tau$-quantile contours} for $Y$, while $\RR_P (Y)$,   $\RR_P (Y)/\|\RR_P (Y)\|$ and $\|\RR_P(Y)\|$ play the roles of vector ranks, signs, and ranks, respectively (\cite{HP:2002,HP:2004,HP:2005,HP:2006,HP:2008}): we call them spherical {\it vector ranks}, { \it signs}, and {\it ranks}.   On the other hand, we call the inverse map $u \longmapsto \QQ_P(u)$ of the vector rank map $y \longmapsto \RR_P (y) = (y/\|y\|)G(\|y\|)$ the {\it vector quantile} map. In both cases, the  relevance of   ranks and signs, whether traditional or spherical, is related to their role as  maximal invariants under   groups  of transformations minimally generating~$\mathcal{P}^1$ or the family~$\mathcal{P}^d_{\mbox{\scriptsize sph}}=\{P_{0,I,f}\}$ of { spherical} distributions, of which distribution-freeness of~$\RR_P$ is just a by-product, as explained in \cite{HW:2003}.
We argue that an adequate depth function, when restricted to those  two particular cases, should lead to the same well-established concepts---classical quantiles, ranks and signs for $\mathcal{P}^1$, and spherical  ones for~$\mathcal P^d_{\mbox{\scriptsize sph}}$---hence 
 should coincide with halfspace depth.

Now, a  closer look at  those two particular cases reveals that  halfspace depth contours, in $\mathcal{P}^1$ and $\mathcal P^d_{\mbox{\scriptsize sph}}$, are the  images, by the vector quantile map $\QQ_P$, of the hyperspheres~$\mathcal{S}(\tau)$ with radii~$\tau\in[0,1)$ centered at the origin.
%
The map $\QQ_P$   is  the
 gradient of a convex function and it transports the spherical uniform distribution $U_d$ on the unit ball $\mathbb S^d$ of~$\R^d$ into the univariate distribution $P~\!\in~\!\mathcal{P}^1$ or into the spherical distribution~$P=P_{0 ,I , f}$ of interest.

 For the case of general distributions $P$, we proceed similarly, and 
  define the map $\QQ_P$
 as  a gradient of a convex function that  transform the spherical uniform distribution~$U_d$ into the target distribution, namely if $U \sim U_d$ then $Y = \QQ_P(U) \sim P$.  It follows by McCann's~\cite{McCann1995} extension of Brenier's celebrated Polar Factorization Theorem~\cite{brenier} that, for any distribution $P$ on~$\R^d$,  such a  gradient $\QQ_P$ exists, and  is   essentially unique. Moreover, when~$P$ has finite moments of order two,  that mapping $\QQ _P$ is the Monge-Kantorovich  optimal transport map that transfers the spherical uniform distribution~$U_d$ to $P$, where optimality is the sense of minimizing  the expected quadratic cost  $\min_\QQ \E_{U} ( \QQ(U) -U)^2$ subject to $U \sim U_d$ and~$\QQ(U) \sim P$.

This suggests a new concept of statistical depth,  which we call the  {\it Monge-Kantorovich (or MK) depth} $\DD^{\scriptsize{\text{MK}}}$, the contours  of which are obtained as the images by $\QQ _P$  of the hyperspheres with radius~$\tau\in[0,1]$. When restricted  to $\mathcal{P}^1$ or $\mathcal P^d_{\mbox{\scriptsize sph}}$, Monge-Kantorovich and halfspace depths coincide. Under suitable regularity conditions due to Caffarelli (see \cite{villani1}, Section~4.2.2),~$\QQ _P$ is a homeomorphism, and its inverse $\RR _P:=\QQ_P^{-1}$  is also the gradient of a convex function; the Monge-Kantorovich depth contours are continuous and the corresponding depth regions are nested, so that Monge-Kantorovich depth indeed provides a center-outward ordering of~$\R^d$, namely,
\begin{equation}\label{MKorder}
y_2 \geq_{\DD_P^{\scriptsize{\text{MK}}}} y_1 \mbox{ if and only if } \Vert\RR_P(y_2)\Vert \leq \Vert\RR _P(y_1)\Vert .
\end{equation}
Thus, our approach based on the theory of measure transportation allows us to define
\begin{enumerate}[(a)] \item an {\em MK vector quantile} map $\QQ_P$, and the associated  {\em MK quantile} correspondence, which maps $\tau\in[0,1]$ to $\QQ_P(\mathcal S(\tau))$, \item an {\em MK vector rank} (or {\em MK signed rank}) function $\RR_P$, which can be decomposed into an {\em MK rank} function $r_P$ from $\R^d$ to $[0,1]$, with $r_P(x):=\|\RR_P(x)\|$, and an {\em MK sign} function~$u_P$, mapping $x\in\R^d$ to $u_P(x):=\RR_P(x)/\|\RR_P(x)\|\in\mathcal S^{d-1}$.
\end{enumerate}
To the best of our knowledge, this is the first proposal of a depth concept based on
the Monge-Kantorovich theory of measure transportation ---hence the first attempt to provide a measure-driven ordering of $\R^d$ based on measure transportation theory.  Previous proposals have been made, however, of measure transportation-based vector quantile functions in Ekeland, Galichon and Henry \cite{EGH} and Galichon and Henry \cite{GH:2012} (with moment conditions) and Carlier, Chernozhukov and Galichon \cite{CCG} (dropping moment conditions) who also extended the notion to vector quantile regression, creating a vector analogue of Koenker and Basset's \cite{KB:78}
scalar quantile regression. More recently, Decurninge \cite{Decurninge}
proposed a new concept of multivariate $L^p$ moments based upon a similar notion. In these contributions, however, the focus is not statistical depth and the associated ranks and quantiles, and the leading case for the reference distribution is uniform on the unit hypercube in $\R^d$, as opposed to the spherical uniform distribution $U_d$ we adopt here as leading case, while pointing out that other reference distributions may be entertained, such as the standard Gaussian distribution on $\R^d$ or the uniform on the hypercube $[0,1]^d$ as mentioned above.

We then proceed to define the empirical notions corresponding to the concepts given above.
We define the empirical MK  vector quantiles and ranks as the essentially unique gradients $\hat\QQ_n$ and $\hat \RR_n$ of a pair of convex functions  solving the Kantorovich dual problem for the Monge optimal transport with quadratic costs.  Using the plug-in principle, we then define the empirical rank and sign maps  as $\|\hat \RR_n\|$ and $\hat \RR_n/\|\hat \RR_n\|$ and the empirical
$\tau$-quantile sets and contours as $\hat \QQ_n(\mathbb{S}(\tau))$ and $\hat \QQ_n(\mathcal{S}(\tau))$.  We establish the uniform convergence of these quantities to their theoretical counterparts.  We derive these results as a consequence of the uniform convergence of empirical transport (vector quantile and rank) maps $\hat \QQ_n$ and $\hat \RR_n$ to their theoretical counterparts $\QQ_P$ and $\RR_P$ on compact subsets of the domain's interior. This is the main theoretical result of the paper presented in Theorem 3.1.
This result in turn is derived through an application of the extended continuous mapping theorem and a set of new theorems on stability of transport under
deterministic perturbations of the source and target measures, given as Theorems A.1 and A.2 in the Appendix, which are new results of independent interest.  Application of the extended continuous mapping theorem allows to us then to replace the deterministic perturbations by stochastic perturbations of measures and obtain the stochastic uniform convergence of the empirical transport maps.

\subsection*{Notation, conventions and preliminaries}

Let $(\Omega, \mathcal{A}, \P)$ be some probability space. Throughout, $\mathcal P$ denotes a class of probability distributions over $\R^d$---unless otherwise specified,  the class of all Borel probability measures on $\R^d$. Denote by $\mathbb S^d:=\{x\in\R^d:\;\|x\|\leq1\}$ the unit ball, and by $\mathcal S^{d-1}:=\{x\in\R^d:\;\|x\|=1\}$ the unit sphere, in $\R^d$. For~$\tau\in(0,1]$, $\mathbb{S}(\tau ):=\{x\in\R^d:\;\|x\|\leq\tau\}$ is the ball, and $\mathcal{S}(\tau ):=\{x\in\R^d:\;\|x\|=\tau\}$ the sphere, of radius $\tau$. Let~$P_X$ stand for the distribution of the random vector $X$. The symbol $\partial$ will denote either the boundary of a set or the subdifferential, as will be clear from the context.
Following Villani \cite{villani1}, we denote by~$g\#\mu$ the {\em image measure} (or {\em push-forward}) of a measure $\mu\in\mathcal P$ by a measurable map~$g:\R^d\rightarrow\R^d$. Explicitly,  for any Borel set $A$, $g\#\mu(A):=\mu(g^{-1}(A))$.
For a Borel subset $\mathbb{D}$ of a vector space equipped with the norm $\|\cdot\|$ and $f: \mathbb{D} \mapsto \R$, let
$$\|f\|_{\mathrm{BL}(\mathbb{D})} := \sup_{x} |f(x)|   \vee \sup_{x \neq x'}  | f(x) - f(x')|\|x - x'\|^{-1} .$$
For two probability distributions $P$ and $P'$ on a measurable space $\mathbb{D}$, define the bounded Lipschitz metric as
$$d_{\BL} (P, P') := \| P - P'\|_{\BL} := \sup_{\|f\|_{\mathrm{BL}(\mathbb{D})} \leq 1} \int f d (P - P'),$$
which metrizes the topology of weak convergence.
Throughout the paper, we let $\mathcal{U}$ and $\mathcal{Y}$ be convex subsets of~$\R^d$ with non-empty interiors. A {\em convex} function $\psi$ on $\mU$ refers to a function
$\psi:\mU \rightarrow\R\cup\{+\infty\}$ for which $\psi((1-t)x+tx')\leq (1-t)\psi(x)+t\psi(x')$ for any $(x,x')$ such that $\psi(x)$ and $\psi(x')$ are finite and for any $t\in(0,1)$. Such a function is continuous on the interior of the convex set  dom $\psi:=\{x\in\mU: \psi(x)<\infty\}$, and differentiable Lebesgue-almost everywhere in dom $\psi$, by Rademacher's theorem.
Write $\nabla\psi$ for the gradient of $\psi$.  For any function $\psi: \mU \mapsto \R \cup \{+\infty\}$, the {\em conjugate} $\psi^*: \mY \mapsto \R \cup \{+\infty\}$ of $\psi$  is defined for each $y \in \mY$ by
$$\psi^*(y) := \sup_{z \in \mU} [y^\top z - \psi(z)].$$ The conjugate~$\psi^\ast$ of $\psi$ is a convex lower-semi-continuous function on $\mY$. We shall call a {\em conjugate pair of potentials} over $(\mU, \mY)$ any pair of lower-semi-continuous convex functions $(\psi, \psi^*)$ that are conjugates of each other.
The transpose of a matrix $A$ is denoted~$A^\top$.
Finally, we call {\em weak order} a complete reflexive and transitive binary relation.
Finally, recall the definition of Hausdorff distance between two non-empty sets $A$ and $B$ in~$\R^d$:
$$
d_H(A,B) :=   \sup_{b \in B}\inf_{a \in A}\| a- b \| \vee \sup_{a \in A}\inf_{b \in B}\| a- b \|.
$$

\subsection*{Outline of the paper}
Section~\ref{sec:MK} introduces and motivates the concepts of statistical depth, vector quantiles and vector ranks based on optimal transport maps. Section~\ref{sec:emp} describes estimators of depth contours, quantiles and ranks, and proves consistency of these estimators. Section~\ref{sec:comp} describes computational characterizations. The appendix presents additional theoretical results and proofs.


\section{Statistical depth and vector ranks and quantiles}
\label{sec:MK}

\subsection{Statistical depth, regions and contours}
The notion of statistical depth serves to define a center-outward ordering of points in the support of a distribution on $\R^d$, for~$d>1$. As such, it emulates the notion of quantile for distributions on the real line. We define it as a real-valued index on $\R^d$ as follows.

\begin{definition*}[Depth and ordering]
\label{def:depth}
A depth function is an upper-semi-continuous mapping  $\DD : \R^d \longmapsto \R$.
In our context these functions will be indexed by a distribution~$P$. The quantity $\DD_P(x)$ is called the {\em depth of} $x$ {\em relative to} $P$. For each $P\in~\!\mathcal P$, the {\em depth ordering} $\geq_{\DD_P}$ {\em associated with} $\DD_P$ is the weak order on $\R^d$ defined, for~$(y_1,y_2)\in~\!\R^{2d}$, by
 \[ y_1\geq_{\DD_P}y_2 \mbox{ if and only if } \DD_P(y_1)\geq \DD_P(y_2),\]
in which case $y_1$ is said to be {\em deeper} than $y_2$ relative to $P$.
\end{definition*}
The depth function thus defined allows graphical representations of the distribution~$P$ through depth contours, which are collections of points of equal depth relative to $P$.
\begin{definition*}[Depth regions and contours]
\label{def:rc}
Let $\DD_P$ be a depth function relative to distribution $P$ on $\R^d$.  The {\em region of depth} $d$  is  the upper contour set of level $d$ of $\DD_P$, namely $\mathbb C_P(d)=\{x\in\R^d:\;\DD_P(x)\geq d\}$; the {\em contour of depth} $d$  is the boundary~$\mathcal C_P(d)= \partial \mathbb C_P(d)$.
\end{definition*}
By construction, the depth regions are nested:  \[\forall (d,d')\in\R_+^2,\;
d'\geq d \;\Longrightarrow\;\mathbb C_P(d')\subseteq\mathbb C_P(d).\] Hence, the depth ordering qualifies as a {\em center-outward ordering} of points in $\R^d$ relative
to the center given by the set of the deepest points,  $\arg\sup_{ x \in \R^d} \DD_P(x).$

It is often convenient to work with depth regions indexed by their probability content.
\begin{definition*}[Depth regions with probability content $\tau$]
\label{def:rcq}For $\tau \in [0,1]$, the depth region with probability content at least $\tau$ is
$$
\mathbb K_P(\tau):=\mathbb{C}_P(d(\tau)),  \quad d(\tau):= \inf\{d \in \R:
\Pr_P (\mathbb{C}(d)) \geq \tau \};
$$
the corresponding contour region is the boundary $\mathcal K_P(\tau):= \partial \mathbb K_P(d)$.

\end{definition*}

\subsection{Liu-Zuo-Serfling axioms and Tukey's halfspace depth}\label{sec:LZS}
\begin{figure}[hhtbp]\label{fig:Tban}
	\begin{center}
	\includegraphics[width=5in, height=5in]{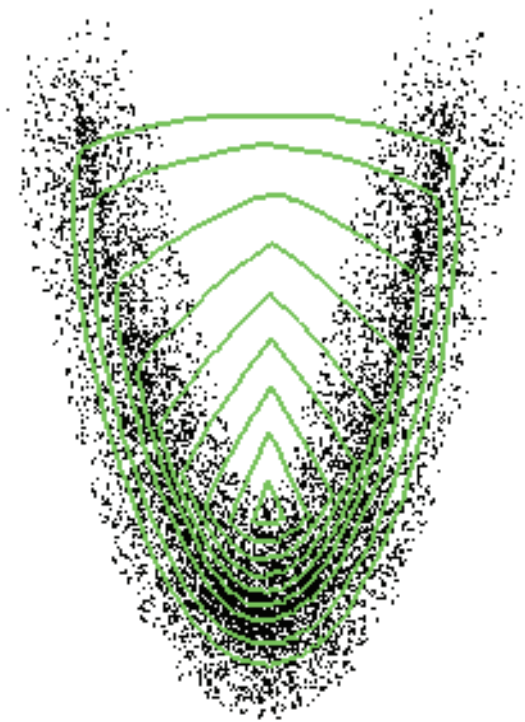}
	\caption{Tukey halfspace depth contours for a banana-shaped distribution, produced with the algorithm of Paindaveine and \v Siman \cite{PS:2012} from a sample of~9999 observations. The banana-like geometry of the data cloud is not picked by the convex contours, and the deepest point is close to the boundary of the support.}
	\end{center}
\end{figure}

The four axioms proposed by Liu \cite{Liu:90} and Zuo and Serfling \cite{ZS:2000} to unify the diverse depth functions proposed in the literature are the following.
\begin{enumerate}
\item[(A1)]  (Affine invariance) $\DD_{ P_{AX+b}}(Ax+b)=\DD_{{ P}_X}(x)$ for any $x\in\mathbb{R}^d$, any
full-rank  $d\times d$ matrix $A$, and any~$b\in\R^d$.
\item[(A2)] (Maximality at the center) If $x_0$ is a center of symmetry for ${ P}$ (symmetry here can be either {\it central}, {\it angular} or {\it halfspace} symmetry), it is {\it deepest}, that is, $\DD_{ P}(x_0)=\max_{x\in{\R^d}}\DD_{ P}(x)$.
\item[(A3)] (Linear monotonicity relative to the deepest points)
If $\DD_{ P}(x_0)=\max_{x\in{\R^d}}\DD_{ P}(x)$, then $\DD_{ P}(x) \leq \DD_{ P}((1-\alpha) x_0 + \alpha x)$ for all $\alpha\in[0,1]$ and $x\in\R^d$: depth is monotonically decreasing  along any straight line running  through a deepest point.
\item[(A4)] (Vanishing at infinity) $\lim_{\Vert x\Vert \to\infty}\DD_{ P}(x) = 0$.
\end{enumerate}

The earliest and most popular depth function is {\em halfspace depth} proposed by Tukey~\cite{Tukey:75}:
\begin{definition*}[Tukey's halfspace depth]
\label{def:half}\!\!\!\! 
The halfspace depth~$\DD^{\scriptsize{\text{Tukey}}}_{ P}(x)$ of a point $x\in~\!\R^d$ with respect to the distribution $ P_X$ of a random vector $X$ on $\R^d$ is defined~as
\[
\DD^{\scriptsize{\text{Tukey}}}_{ P_X}(x):=\min_{\phi\in\mathcal S^{d-1}} \P [(X-x)^\top \phi\geq 0].
\]
\end{definition*}
Halfspace depth relative to any distribution with nonvanishing density on $\R^d$ satisfies~(A1)-(A4).
The appealing properties of halfspace depth are well known and well documented: see Donoho and Gasko \cite{DG:92}, Mosler \cite{Mosler:2002}, Koshevoy \cite{Koshevoy:2002}, Ghosh and Chaudhuri \cite{GC:2005}, Cuestas-Albertos and Nieto-Reyes \cite{CN:2008}, Hassairi and Regaieg \cite{HR:2008}, to cite only a few.  Halfspace depth takes values in $[0, 1/2]$, and  its  contours  are continuous and convex; the corresponding regions are closed, convex, and nested as $d$ decreases. Under very mild  conditions, halfspace depth moreover fully characterizes the distribution~$P$. For somewhat less satisfactory features, however, see Dutta et al. \cite{DGC:2011}. An important feature of halfspace depth is the convexity of its contours, which implies that halfspace depth contours cannot pick non convex features in the geometry of the underlying distribution, as illustrated in Figure~1.

We shall propose below a new depth concept, the Monge-Kantorovich (MK) depth, that relinquishes the affine equivariance and star convexity of contours imposed by Axioms~(A1) and~(A3) and recovers non convex features of the underlying distribution.
As a preview of the concept, without going through any definition, we illustrate in Figure 2 (using the same banana-shaped distribution as in Figure 1) the ability of the
MK depth to capture
non-convexities.  In what follows, we characterize these abilities more formally.  We shall emphasize that this notion
comes in a package with new, interesting notions of vector ranks and quantiles, based on optimal
transport, which reduce to classical notions in the univariate and multivariate spherical cases.
\subsection{Monge-Kantorovich depth}
The principle behind the notion of depth we define here is to map the depth regions and contours relative to a well-chosen reference distribution $F$, into depth contours and regions relative to a distribution of interest $P$ on $\R^d$,
using a well-chosen mapping. The mapping proposed here is
  the gradient
of a convex function $\nabla \psi$  such that if $U$ has distribution $F$, then $Y = \nabla \psi (U)$ has distribution $P$, or, in terms of measures, $\nabla \psi \# F = P$.
The gradient $\nabla \psi$ is said to to push $F$ forward to~$P$, which is conventionally denoted by the push-forward notation, $\nabla \psi \# F = P$, which is defined
in the notation section.  

 The gradient of a convex function property is a generalization of monotonicity in the one-dimensional case. When $F$ and $P$ have finite second-order moments, these maps are the optimal Monge-Kantorovich transport maps from $F$ to $P$ for the quadratic cost, as explained below.
 In the unidimensional case, when $F$ is the standard uniform, the gradient/optimal transport map $\nabla \psi$ coincides with the classical quantile function.
\begin{figure}[htbp]\label{fig:MKban}\begin{center} \includegraphics[width=4in, height=5in]{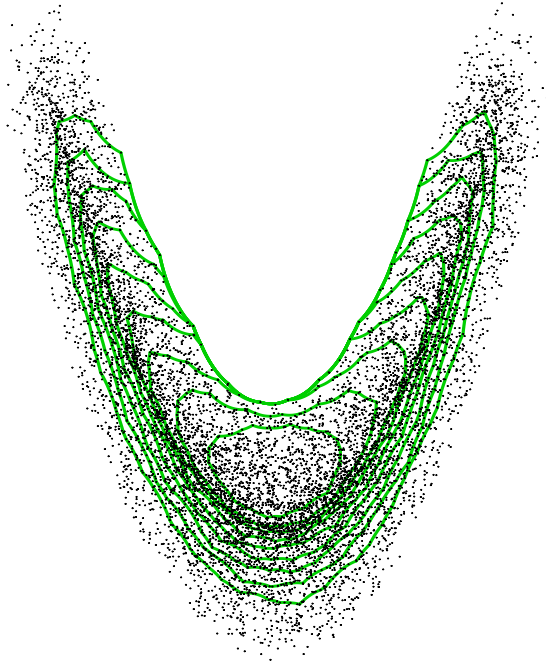}\caption{The Monge-Kantorovich depth contours for the same banana-shaped distribution from a sample of~9999 observations, as in Figure~1. The banana-like geometry of the data cloud is correctly picked up by the non convex contours.}
	\end{center}
\end{figure}

The following theorem, due to Brenier \cite{brenier} and McCann \cite{McCann1995}, establishes existence of gradients of convex functions with the required properties.
\begin{theorem}[Brenier-McCann's Existence Result]
\label{thm:polar}
Let $P$ and $F$ be two distributions on $\R^d$.  (1) If $F$ is absolutely continuous with respect to the Lebesgue measure on~$\Bbb{R}^d$, with support contained in a convex set $\mU$, the following holds:  there exists a convex function $\psi:\mU\rightarrow\R\cup\{+\infty\}$ such that $\nabla\psi\# F=P$.  The function $\nabla\psi$ exists and is unique, $F$-almost everywhere. (2) If, in addition, $P$ is absolutely continuous on $\Bbb{R}^d$ with support contained in a convex set $\mY$, the following holds: there exists a convex function~$\psi^*:\mY\rightarrow\R\cup\{+\infty\}$ such that $\nabla\psi^*\# P=F$.  The function $\nabla\psi^*$ exists, is unique and equal to $\nabla \psi^{-1}$, $P$-almost everywhere. \end{theorem}

\begin{remark}[Interpretation as a Monge-Brenier Optimal Transport]   If $P$ and $F$ have finite  second moments,
$\QQ$ is  $F$-almost everywhere equal to  the {\em optimal transport plan}  $\nabla \psi$  from $F$ to $P$ for quadratic cost:  namely, the map  $\QQ: \R^d\longrightarrow\R^d$   solves the problem 
\begin{eqnarray*}
\inf_{Q}  \int (u - Q(u))^2 d F(u) :  \quad Q\#F=P,
\end{eqnarray*}
or, equivalently,
\begin{eqnarray}
\sup_{Q} \int u^\top Q(u)\;dF(u) :  \quad Q\#F=P.\label{eq:OT}
\end{eqnarray}
This definition has a classical counterpart in the case of univariate distributions.
 When $d=1$ and $F$ is uniform on $[0,1]$, the optimal transport $u\mapsto \QQ(u)$ is the classical quantile function for distribution~$P$. \qed
 \end{remark}

We now state a fundamental duality result due to Kantorovich and Brenier, which we explicitly rely on in Section~\ref{sec:emp}.

\begin{theorem}[Kantorovich-Brenier, see \cite{villani1}]\label{dual}
Suppose hypothesis (1) of Theorem~\ref{thm:polar} holds and  $P$ and $F$ have finite  second moments, then the function $\psi$, or optimal potential,  solves the optimization problem
\begin{eqnarray}
\label{eq:Kdual}
\int\psi dF+\int\psi^\ast dP=\inf_{(\varphi, \varphi^*)} \left( \int\varphi dF+\int\varphi^\ast dP \right),
\end{eqnarray}
where the infimum is taken over the class of  conjugate pairs of potentials $(\varphi, \varphi^*)$
over~$( \mU, \mY)$.
\end{theorem}

\begin{remark} This problem is dual to the optimal transport problem (\ref{eq:OT}). Moreover,  under the hypotheses of Theorem~\ref{dual}, $\nabla\psi$ is the unique  optimal transport map from~$F$ to~$P$ for quadratic cost, in the sense that any other optimal transport coincides with $\nabla \psi$ on a set of $F$-measure 1 (see \cite{villani1}). Under the hypotheses of Theorem~\ref{dual} and hypothesis~(2) of Theorem~\ref{thm:polar},  $\nabla\psi^*$ is the unique optimal (reverse) transport map from $P$ to~$F$ for quadratic cost, in the sense that any other optimal transport coincides with $\nabla \psi^*$ on a set of $P$-measure one (see \cite{villani1}). \qed
\end{remark}

Next we use Theorem~\ref{thm:polar} to define a natural notion of \textit{vector quantiles} and \textit{vector ranks}.

\begin{definition}[Monge-Kantorovich vector quantiles and ranks]
\label{def:MK}
Let $F$ be an absolutely continuous reference distribution with support in a convex set $\mU \subseteq \R^d$, and let~$P$ be an arbitrary distribution with support in a convex set $\mY \subseteq \R^d$.   Let $\nabla \psi$ be the $F$-almost surely unique gradient of a convex function $\psi$
of Theorem~\ref{thm:polar} and let~$\psi^\ast$ be the conjugate of $\psi$ over $(\mU, \mY)$. Vector quantiles and ranks  are defined as follows:
\begin{equation*}\QQ_P(u) \in \arg \sup_{y \in \mY}[ y^\top u - \psi^*(y)] ,  \  \ { u \in \mU }; \quad
\RR_P(y) \in \arg \sup_{u \in \mU} [y^\top u - \psi(u)], \ \ { y \in \R^d}.
\end{equation*}
\end{definition}
\begin{remark} Thus we define the MK vector quantiles $\QQ_P$ and ranks $\RR_P$ as any solutions of the optimization problems in the display above.  Our definition here does not impose any moment condition and ensures that the quantities are defined for every value of the argument in the appropriate domains. By the envelope theorem and Rademacher's theorem (\cite{villani1}),  the  maps $\QQ_P$ and $\RR_P$ essentially coincide with the gradients $\nabla \psi$ and $\nabla \psi^*$ of conjugate potentials $\psi$ and $\psi^*$, namely
\begin{equation}
\QQ_P = \nabla \psi \text{ a.e.  on } \ \mathcal{U}, \quad \RR_P = \nabla \psi^* \text{ a.e.  on } \ \mathcal{Y},
 \end{equation}
where ``a.e." abbreviates ``almost everywhere with respect to the Lebesgue measure". In the fact, the equality
holds everywhere on certain domains under condition (C) stated below.  Under the conditions of Theorem~\ref{dual}, the pair $(\psi, \psi^*)$ has the variational characterization given in~(\ref{eq:Kdual}). \qed\end{remark}

When requiring regularity of vector quantiles and ranks, we shall impose the following condition on the conjugate pair of optimal potentials $(\psi,\psi^\ast)$ over $(\mU, \mY)$.
\begin{itemize}
\item[(C)]
Let $\mathcal{U}$ and $\mathcal{Y}$ be closed, convex subsets of $\R^d$, and $\mU_0 \subset \mU$ and $\mY_0 \subset \mY$ 
be open, non-empty sets in $\R^d$.   Let $\psi: \mathcal{U} \mapsto \R$ and $\psi^*: \mathcal{Y} \mapsto \R$ form a conjugate pair over~$(\mU, \mY)$ and possess gradients $\nabla \psi(u)$ for all $u \in \mU_0$, and $\nabla \psi^*(y)$ for all $y \in \mY_0$. The gradients $\nabla \psi |_{\mU_0}:  \mU_0 \mapsto \mY_0$ and $\nabla \psi^* |_{\mY_0}: \mY_0 \mapsto \mU_0$ are homeomorphisms and~$\nabla \psi|_{\mU_0} = (\nabla \psi^* |_{\mY_0})^{-1}$.
\end{itemize}
Under Condition (C), we have: \begin{equation}\label{as gradients}
\QQ_P(u) = \nabla\psi(u)  \text{ for all } u \in \mU_0, \quad
\RR_P(y) = \nabla\psi^\ast(y)=\left( \nabla\psi \right)^{-1}(y) \text{ for all } y \in \mY_0,
\end{equation}
that is,  vector ranks and quantiles are defined as gradients of conjugate potentials for each (as opposed to almost every) value in the indicated sets, and inverse functions of each other.

Sufficient conditions for Condition (C) in the context of Definition~\ref{def:MK} are provided by Caffarelli's regularity theory (Villani \cite{villani1}, Theorem 4.14). One set of sufficient conditions is as follows.
\begin{lemma}[Caffarelli's Regularity, \cite{villani1}, Theorem 4.14]
\label{prop:C}
Suppose that $P$ and $F$ admit densities,  which are of smoothness class $C^{\beta}$ for $\beta>0$ on convex, compact support sets ${\rm cl}(\mY_0)$ and ${\rm cl}(\mU_0)$,
and the densities are bounded away from zero and above uniformly on the support sets. Then Condition (C) is satisfied for the conjugate pair $(\psi,\psi^\ast)$ such that~$\nabla\psi\#F=P$ and $\nabla\psi^\ast\#P=F.$
\end{lemma}

We now can give our {\em main} definition -- that of multivariate notions of quantiles and ranks, through which a depth function will be inherited from the reference distribution~$F=U_d$.
\begin{definition}[\textbf{Monge-Kantorovich depth, quantiles, ranks and signs}]
\label{def:MK} Let~$F$ be the spherical uniform distribution $U_d$ on a unit ball $\mU = \mathbb{S}^d$, and $P$ be an arbitrary distribution with support in  a convex region $\mY \subseteq \R^d$. MK quantiles, ranks, signs and depth are defined as follows.
\begin{enumerate}
\item  The {\em MK rank} of $y \in \R^d$ is $\|\RR_P(y)\|$ and the {\em MK sign} is~$\RR_P(y)/\|\RR_P(y)\|$.

\item   The {\em MK} $\tau$-{\em quantile contour} is the set $\QQ_P(\mathcal S(\tau))$ and the {\em MK depth region} with probability content $\tau$ is $\QQ_P(\mathbb S(\tau))$.

\item The {\em MK  depth} of $y \in\R^d$ is  the depth of $\RR_P(y)$ under $\DD^{\scriptsize{\text{Tukey}}}_{U_d}$:
\[
\DD^{\rm\scriptsize MK}_P(y):=\DD^{\scriptsize{\text{Tukey}}}_{U_d}(\RR_P(y)).
\]
\end{enumerate}
\end{definition}
The notion of depth proposed in Definition~\ref{def:MK} is based on an optimal transport map from the reference spherical uniform distribution $F=U_d$ to the distribution of interest~$P$.  
Under Condition (C), $\QQ_P$ and $\RR_P$ are continuous and are mutual inverse maps, so that the MK $\tau$-quantile contours are continuously deformable into spheres and the MK depth regions with probability content $\tau$ are nested.

By choosing other reference distributions $F$, such as the uniform distribution on a unit hypercube,  or the standard Gaussian distribution,
 we can give a more general definition of MK ranks, quantiles, and signs, which may be of interest.

\begin{definition}[Monge-Kantorovich depth, quantiles, ranks and signs for general~$F$]
\label{def:MK2}Let $F$ be an absolutely continuous reference distribution with support contained in a convex region $\mU \subseteq \Bbb{R}^d$, and let $\| \cdot \|$ be a norm on $\mU$. Let $\DD_F: \R^d \to \Bbb{R}_+$ be an associated reference depth function and $\mK(\tau)$ the associated $\tau$-quantile contour and $\bK(\tau)$
the associated depth region with probability content $\tau$.  The MK quantiles, ranks, signs and depth are defined as follows.
\begin{enumerate}
\item  The {\em MK rank} of $y \in \R^d$ is $\|\RR_P(y)\|$ and the {\em MK sign} is~$\RR_P(y)/\|\RR_P(y)\|$.
\item   The {\em MK} $\tau$-{\em quantile} is the set $\QQ_P(\mK(\tau))$ and the {\em MK depth
region} with probability mass $\tau$ is $\QQ_P(\bK(\tau))$.
\item The {\em MK  depth} of $y \in\R^d$ is the depth of $\RR_P(y)$ under $\DD_F$:
$$
\DD^{\rm\scriptsize MK}_P(y):=\DD_F(\RR_P(y)).
$$
\end{enumerate}
\end{definition}
Of course, all the quantities thus defined depend on the choice of the reference distribution $F$ and the depth function $\DD_F$. 
\begin{remark}  When the reference distribution~$F$ is spherical, it is natural to use Tukey's depth function~$\DD_F = \DD^{\scriptsize{\text{Tukey}}}_{F}$  to define the MK depth of $y \in\R^d$ relative to~$P$ as the halfspace depth of $\RR_P(y)$ relative to the reference distribution~$F$, namely $$\DD^{\rm\scriptsize MK}_P(y):=\DD^{\scriptsize{\text{Tukey}}}_{F}(\RR_P(y)).$$
The choice of halfspace depth may be less natural for non-spherical reference distributions. One example is where $F$ is the standard uniform distribution $U[0,1]^d$ on the unit cube $[0,1]^d$. Then it seems natural to use the sup norm $\|\cdot\|_\infty$ as the norm $\| \cdot \|$ and the depth function $\DD_{U[0,1]^d}(y) = 1/2 - \| y - \mathbf{1}/2\|_\infty$, where $\mathbf{1} = (1,\dots,1)'$, in which case $\bK(\tau)$ is a  cube of diameter $\tau^{1/d}$ centered at $\mathbf{1}/2$. In this
case, the MK depth is
$$\DD^{\rm\scriptsize MK}_P(y):=\DD_{U[0,1]^d}(\RR_P(y)). $$
 \end{remark}

\subsection{Monge-Kantorovich depth with spherical uniform reference distribution}

Here we consider in more detail the Monge-Kantorovich depth defined from a baseline spherical uniform distribution $U_d$ supported on the unit ball $\mathbb S^d$ of $\R^d$. Recall that this distribution is
that  of a random vector $r\phi$, where $r$ is uniform on $[0,1]$, $\phi$ is uniform on the unit sphere~$\mathcal S^{d-1}$, and $r$ and $\phi$ are mutually independent.

The spherical symmetry of distribution~$U_d$ produces halfspace depth contours that are concentric spheres, the deepest point being the origin.
The radius $\tau $ of the ball $\mathbb{S}(\tau )=\{x\in\R^d:\;\|x\|\leq\tau\}$ is also its $U_d$-probability contents, that is,~$\tau =U_d({{\mathbb{S}}}(\tau ))$.
Letting $\theta :=\arccos \tau $, the halfspace depth with respect to $U_d$  of  a point~$\tau u\in\mathcal{S}(\tau ):=\{x\in\R^d:\;\|x\|=\tau\}$, where~$\tau \in (0,1]$ and $u\in\mathbb S^d$, is
\begin{equation}\label{hlfd} \DD_U(\tau u)=\left\{
   \begin{array}{ll}   \pi^{-1}[\theta -\cos\theta \log \vert \sec\theta + \tan\theta\vert ] &
   d\geq 2\\(1-\tau )/2 & d=1.\end{array}\right.
\end{equation}
Note that for $d=1$, $u$ takes values $\pm 1$ and, in agreement with rotational symmetry of~$U_d$, that depth does not depend on $u$.

The principle behind the notion of depth we investigate further here is to map the depth regions and contours
relative to the  spherical uniform distribution $U_d$, namely, the concentric spheres, into depth contours and regions relative to a distribution of interest~$P$ on~$\R^d$
using the optimal transport plan from~$U_d$ to~$P$.
Under the sufficient conditions for Condition (C) provided in Lemma~\ref{prop:C} (note that the conditions on $F$ are automatically satisfied in case $F=U_d$), $\QQ_P$ and $\RR_P$ are continuous and are inverse maps of each other, so that the MK depth contours are continuously deformable into spheres, the MK depth regions are nested, and regions and contours, when indexed by probability content, take the  respective forms
\[
\QQ_P\left( \mathbb S(\tau)\right) \mbox{ and } \QQ_P\left( \mathcal S(\tau)\right),\mbox{ for }\tau\in(0,1].
\]

\subsubsection*{MK depth is halfspace depth in dimension 1}
The halfspace depth of a point $x\in\R$ relative to a distribution $P$ over $\R$ takes the very simple form
\[
D^{\scriptsize{\text{Tukey}}}_P(x)=\min(P(x), 1-P(x)),
\]    where, by abuse of notation, $P$ stands for both distribution and distribution function.
The nondecreasing map defined for each $x\in\R$ by $x\mapsto \RR_P(x)=2P(x) -1$ is the derivative of a convex function and it transports distribution $P$ to $U_1$, which is uniform on~$[-1,1]$, i.e., $\RR_P\#P=U_1$. Hence $\RR_P$ coincides with the MK vector rank of Definition~\ref{def:MK}. Therefore, for each $x\in\R$,
\[
D_P(x)=D^{\scriptsize{\text{Tukey}}}_{U_d}(\RR_P(x))=\min(P(x), 1-P(x))
\]
and MK depth coincides with Tukey depth in case of all distributions with nonvanishing densities on the real line.

More generally (still in the univariate case), denoting by $F_1$ and $F_2$ the distribution functions associated with two absolutely continuous distributions $P_1$ and~$P_2$, the mapping $F_2^{-1}\circ F_1$, being  monotone increasing,  is also the optimal transport from $P_1$ to~$P_2$. The same transformation has been studied, in a different context, by Doksum~\cite{Doksum:74} and Doksum and Sievers~\cite{DS:76}; see also the concept of  convex ordering proposed by van Zwet~\cite{VZ:1964}.

\subsubsection*{MK depth is halfspace depth for elliptical families}
As explained in the introduction, a  $d$-dimensional random vector $X$ has elliptical distribution $P_{\mu,\Sigma,g}$ with location~$\mu\in~\!\R^d$, positive definite symmetric $d\times d$ scatter matrix $\Sigma$ and radial density function $g$ (radial distribution function $G$) if  and only if, denoting by $\Sigma^{1/2}$ the symmetric root of $\Sigma$,  \linebreak  $Y:=\Sigma^{-1/2}(X-\mu)$ has spherical distribution $P= P_{0,I,g}$ (hence~$\Vert Y\Vert $ has density $f$), which holds if and only if
\begin{equation}
\label{eq:ell}
\RR_P (Y):=\frac{Y}{\Vert Y\Vert }G\big(\Vert Y\Vert\big)
\mbox{ is distributed according to } U_d.
\end{equation}
Let $\Psi(t) = \int_{-\infty}^t G(r) d r$, and note that the map $z\mapsto \RR_P(z)$ is the gradient of~$\psi^\ast(z):= \Psi(\Vert z \Vert)$ so that, from (\ref{eq:ell}), $\nabla\psi^\ast\#P=U_d$ as Definition~\ref{def:MK} requires.  That  $\psi^\ast$ is convex follows from Theorem 5.1 of \cite{RockafellarConvex} by noting that $\psi^\ast$ is a composition of $\Psi: \R \to \R$, a convex, non-decreasing map, and $\|\cdot \|: \R^d \to R$, a convex function by definition of the norm.   As a consequence, the mapping $\RR_P$ in (\ref{eq:ell}) is  the MK vector rank function associated with $P=P_{0,I,f}$; and,   the MK depth contours (with probability content $\tau$) of~$P$ are spheres with radii $G^{-1}(\tau)$ centered at the origin:
$$
D_P(x) = \{ y \in \R^d:  \| y \| \leq G^{-1} (\tau) \}.
$$
These spheres are halfspace depth contours for $P$.  This is the precise sense in which MK depth reduces to halfspace for elliptical families.

It should be noted above, that we treat location and scatter parameters as known, and transform $X$ to a vector $Y$ in isotropic position.  This transformation
ensures basic invariance properties of the resulting depth, ranks, and quantiles with respect to affine transformations. When those parameters are unknown, they will have to be replaced with by affine-equivariant estimators, as in the usual definition of elliptical ranks and signs (see,e.g., \cite{HP:2002})  in order 
to insure similar invariance properties for the empirical analogs.  Without the aforementioned transformation, however, the invariance properties are not guaranteed, owing to the fact that composition of two gradients of convex functions is not necessarily the gradient of a convex function, unless the composition has a specific structure, as is the case above.


\section{Empirical depth, ranks and quantiles}
\label{sec:emp}

Having defined Monge-Kantorovich vector quantiles, ranks and depth relative to a distribution $P$ based on reference distribution $F$ on $\R^d$, we now turn to the estimation of these quantities. Hereafter, we shall assume that Condition (C) holds. Then, the MK vector quantiles and ranks of Definition~\ref{def:MK} are
\begin{equation}
\QQ_P(u) := \nabla\psi(u) , \quad
\RR_P(y) := \nabla\psi^\ast(y)=\left( \nabla\psi \right)^{-1}(y),
\end{equation}
for each $u \in \mU_0$ and $y \in \mY_0$, respectively. We define $\Phi_0(\mU, \mY)$ as a collection of conjugate potentials $(\phi, \phi^*)$
on $(\mU, \mY)$ such  that $\phi(u_0)= 0$ for some fixed point $u_0 \in  \mU_0$. Under
the conditions of Theorem~\ref{dual},  the potentials $(\psi, \psi^*)$ solve
the dual problem \begin{equation}\label{kantorovich 2}
\int \psi dF+ \int \psi^* d P = \inf_{ (\varphi, \varphi^*) \in \Phi_0(\mU, \mY) } \int \varphi dF + \int \varphi^* dP .
\end{equation}
Constraining the conjugate pair to lie in $\Phi_0(\mU, \mY)$ is a normalization that (without any loss of generality) pins down the constant, so that $(\psi,\psi^\ast)$ are uniquely determined, as argued in the proof.

We propose empirical versions of MK quantiles and ranks based on estimators $\hat P$ of~$P$.  The typical case is when the reference measure $F$ is known.  However, our theory  allows us to handle the case where $F$ is  itself unknown, and so it is estimated by some $\hat F$.    This is indeed useful for at least two reasons.  First,
we may be interested in a classical problem of comparing one distribution $P$ to a reference distribution $F$, both of which are known only up to a random sample available
from each of them. Second, we may be interested in discretizing $F$ for computational reasons, as we discuss in Section~\ref{sec:comp}, in which case the discretized $F$ is the estimator of $F$.

\subsection{Conditions on estimators of $P$ and $F$}
Suppose
that $\{\hat P_n\}_{n=1}^\infty$  and $\{\hat F_n\}^\infty_{n=1}$  are  sequences of random measures
on $\mY$ and  $\mU$, with finite total mass, that are consistent for $P$ and $F$, in the sense that
\begin{equation}\label{empirical KR}
d_{\KR}(\hat P_n, P) \to_{\P^*} 0, \quad d_{\KR}(\hat F_n, F) \to_{\P^*} 0, \end{equation}
where $\to_{\P^*}$ denotes convergence in (outer) probability under probability measure~$\P$,
see van der Vaart and Wellner \cite{VvW}. A basic example is where $\hat P_n$ is the empirical distribution of a random sample
$(Y_i)_{i=1}^n$ drawn from $P$ and $\hat F_n$ is the empirical distribution of a random sample
$(U_i)_{i=1}^n$ drawn from $F$.  Other, much more complicated examples, including
smoothed empirical measures and data originating from dependent processes, satisfy sufficient conditions for (\ref{empirical KR}) that we now give.
In order to develop some examples,
we introduce an ergodicity condition:
\begin{itemize}
\item[(E)] Let $\mathcal{W}$ be a measurable subset of $\R^d$. A data stream $\{(W_{t,n})_{t=1}^n\}_{n=1}^\infty$, with $W_{t,n} \in \mathcal{W} \subseteq \R^d$ for each $t$ and $n$, is ergodic for the probability law $P_W$ on~$\mathcal{W}$ if for each~$g: \mathcal{W} \mapsto \R$ such that $\|g\|_{\BL(\mathcal{W})} < \infty$, the law of large numbers holds:
\end{itemize}
\begin{equation}\label{erg}
\frac{1}{n} \sum_{t=1}^n g(W_{t,n}) \to_{\P}  \int g(w) d P_W(w).
\end{equation}

The class of ergodic processes is extremely rich, including in particular the following cases:
\begin{itemize}
\item[(E.1)] $W_{t,n}=W_t$, where $(W_t)_{t=1}^\infty$ are independent, identically distributed random vectors
with distribution $P_W$;
\item[(E.2)] $W_{t,n}=W_t$, where $(W_t)_{t=1}^\infty$ is stationary strongly mixing process with marginal distribution $P_W$;
\item[(E.3)] $W_{t,n}=W_{t}$,  where $(W_t)_{t=1}^\infty$ is an  irreducible and aperiodic Markov chain  with invariant distribution $P_W$;
\item[(E.4)] $W_{t,n}=w_{t,n}$, where $(w_{t,n})_{t=1}^n$ is a  deterministic sequence of points such that  (\ref{erg}) holds deterministically.
\end{itemize}
For a detailed motivation and discussion of the use of deterministic sequences such as, for example, the so-called {\it low-discrepancy sequences}:  see, e.g.,  Chapter 9 and, more particularly, page 314 of \cite{ken:judd}.

Thus, if  the data stream $\{(W_{t,n})_{t=1}^n\}_{n=1}^\infty$ is ergodic
for $P_W$, we can estimate $P_W$ by the empirical and smoothed empirical
measures
$$
\hat P_W(A) = \frac{1}{n} \sum_{t=1}^n 1\{W_{t,n} \in A \}, \ \quad
\ \tilde P_W(A) =   \frac{1}{n} \sum_{t=1}^n \int_{\R^d} 1\{W_{t,n} + h_n \eps \in A \cap \mathcal{W} \} d \Phi (\eps),
$$
where $\Phi $ is the probability law of the standard $d$-dimensional Gaussian vector, $N(0,I_d)$,
and  $h_n \geq 0$  a semi-positive-definite matrix of bandwidths such that
$\|h_n\| \to 0$ as~$n \to \infty$.   Note that $\tilde P_W$ may not integrate to 1, since we are forcing
it to have support in $\mathcal{W}$.

\begin{lemma}\label{ergodic} Suppose that $P_W$ is absolutely continuous with support
in the compact set $\mathcal{W} \subset 
\Bbb{R}^d$. If  $\{(W_{t,n})_{t=1}^n\}_{n=1}^\infty$ is ergodic for $P_W$ on
$\mathcal{W}$, then
$$
d_{\BL} (\hat P_W, P_W) \to_{\Pr^*} 0, \quad d_{\BL} (\tilde P_W, P_W) \to_{\Pr^*} 0.
$$
Thus, if $P_Y:=P$ and $P_U:=F$ are absolutely continuous with support sets contained in compact sets $\mY$ and $\mU$,
 and if $\{(Y_{t,n})_{t=1}^n\}_{n=1}^\infty$ is ergodic for $P_Y$ on $\mY$ and~$\{(U_{t,n})_{t=1}^n\}_{n=1}^\infty$
is ergodic for $P_U$ on $\mU$, then  $\hat P_n= \hat P_W$ or $\tilde P_W$ and $\hat F_n = \hat P_U$ or~$\tilde P_U$ obey condition
(\ref{empirical KR}).
\end{lemma}

Absolute continuity of $P_W$ in Lemma~\ref{ergodic} is invoked to show that the smoothed estimator $\tilde P_W$ is asymptotically
non-defective.

\subsection{Empirical vector quantiles and ranks}

We base empirical versions of MK quantiles, ranks and depth on estimators $\hat P_n$ for $P$ and $\hat F_n$ for $F$ satisfying (\ref{empirical KR}). This includes cases where
the reference measure $F$ is known, i.e. $\hat F_n=F$. Recall Assumption~(C) is maintained throughout this section.
\begin{definition}[Empirical Monge-Kantorovich vector quantiles and ranks]
Empirical vector quantile $\hat\QQ_n$ and vector rank $\hat\RR_n$ are any pair of functions satisfying, for each~$u \in \mU$ and $y \in \mY$,
\begin{equation}\label{eq: empirical Q and R}
\hat \QQ_n(u) \in \arg \sup_{y \in \mY}[ y^\top u - \hat \psi_n^*(y)] , \quad
\hat \RR_n(y) \in \arg \sup_{u \in \mU} [y^\top u - \hat \psi_n(u)],
\end{equation}
where  $(\hat \psi_n, \hat \psi^*_n) \in \Phi_0(\mU, \mY)$  is such that
\begin{equation}\label{empirical kantorovich 2}
\int \hat \psi_n d\hat F_n + \int \hat \psi_n^* d \hat P_n =\inf_{ (\phi, \phi^*) \in \Phi_0(\mU, \mY) } \int \varphi d\hat F_n + \int \varphi^* d\hat P_n.
\end{equation}
\end{definition}

We now state the main result of Section \ref{sec:emp}.

\begin{theorem}[Uniform Convergence of Empirical Transport Maps]\label{theorem: empirical}
Suppose that the sets $\mU$ and $\mY$ are compact subsets of $\R^d$, and that the probability measures $P$ and~$F$ are absolutely continuous with respect to the Lebesgue measure, with ${\rm support }(P) \subseteq \mY$ and ${\rm support }(F) \subseteq \mU$.  Suppose
that $\{\hat P_n\}$  and $\{\hat F_n\}$  are  sequences of random measures
on~$\mY$ and  $\mU$, with finite total mass, that are consistent for $P$ and $F$ in the sense
of~(\ref{empirical KR}). Suppose that Condition (C) holds
for the solution of (\ref{kantorovich 2}) for $\mY_0 := {\rm int}({\rm support}(P))$ and~$\mU_0:= {\rm int}({\rm support}(F))$. Then,  as $n \to \infty$, for any closed set $K \subset  \  \mU_0$ and any closed set $K' \subset   \mY_0$,
$$
\sup_{ u \in K} \|  \hat \QQ_n(u) - \QQ_P(u) \|  \to_{\P^*} 0,   \quad
\sup_{ y \in K'} \|  \hat \RR_n(y) - \RR_P(y) \|  \to_{\P^*} 0,
$$ and $$
\sup_{A \subseteq K }d_H(\hat \QQ_{n}(A),\QQ_P(A))  \to_{\P^*} 0, \quad
\sup_{A' \subseteq K'} d_H(\hat \RR_{n}(A'), \RR_P(A'))  \to_{\P^*} 0,
$$
where the suprema are taken over nonempty subsets.
\end{theorem}

The first result establishes
the uniform consistency of empirical vector quantile and rank maps, hence also of empirical ranks and signs. The set $\QQ_P(K)$ with $K = \mathbb{K}(\tau)$ is the statistical depth contour with probability content $\tau$.  The second result, therefore, establishes consistency of the approximation $\hat\QQ_n(K)$ to the theoretical depth region~$\QQ_P(K)$.

\subsection{Empirical MK quantiles, ranks, and signs and their convergence}

We work with the conditions of the previous theorem, but here, for the sake of simplicity,
we first consider the lead case where $F$ is known, i.e. $\hat F_n = F$.

\begin{definition}[Empirical MK depth, quantiles, ranks and signs for known $F$]
\label{def:EMK1}Let~$F$ be an absolutely continuous reference distribution with support contained in a convex region $\mU \subseteq \Bbb{R}^d$, and let $\| \cdot \|$ be a norm on $\mU$. The MK empirical quantiles, ranks, signs and depth are defined as follows.
\begin{enumerate}
\item  The {\em MK empirical rank} and {\em sign} of $y \in \R^d$ are $\|\hat \RR_n(y)\|$ and~$ \hat \RR_n (y)/\|\hat \RR_P(y)\|$.
\item   The {\em MK empirical} $\tau$-{\em quantile contour} is the set $\hat \QQ_n(\mK(\tau))$ and the {\em MK empirical depth region} with probability mass $\tau$ is $\hat \QQ_n(\bK(\tau))$.
\item The {\em MK empirical depth} of $y \in\R^d$ is the depth of $\hat \RR_n(y)$ under $\DD_F$:
$$
\hat\DD^{\rm\scriptsize MK}_{P,n}(y):=\DD_F(\hat \RR_n(y)).
$$
\end{enumerate}
\end{definition}

Uniform convergence of empirical MK  rank, signs and depth to their theoretical counterparts follows by an application of the Extended Continuous Mapping Theorem.
\begin{corollary}
\label{cor:ERS}  Work with the assumptions of Theorem~\ref{theorem: empirical}, and assume that $D_F$ is continuous on $\mU_0$.  As $n \to \infty$, for any closed set $K' \subset   \mY_0$,
\begin{eqnarray*}
\sup_{ y \in K'} |  \|\hat \RR_n(y)\| - \|\RR_P(y)\| |  \to_{\P^*} 0,  \\
\sup_{ y \in K'} |   \hat \RR_n (y)/\|\hat \RR_n(y)\| -  \RR_n (y)/\| \RR_P(y) \| \big |
\to_{\P^*} 0,  \\
\sup_{ y \in K'} |  \hat\DD^{\rm\scriptsize MK}_{P,n}(y) -  \DD^{\rm\scriptsize MK}_{P}(y) |  \to_{\P^*} 0.  \\
 \end{eqnarray*}
\end{corollary}

Uniform convergence of MK empirical $\tau$-quantile contours and MK empirical depth regions with probability content $\tau$ follows also through an application of the Extended Continuous Mapping Theorem.

\begin{corollary} \label{cor:eqd} Work with the assumptions of Theorem~\ref{theorem: empirical}.
Consider $\mathcal{T} \subset (0,1)$ such that $\mathrm{cl}(\cup_{\tau \in \mathcal{T}} \mathbb{K}(\tau)) \subset \mU_0$, 
then
\begin{eqnarray*}\quad \sup_{\tau\in \mathcal{T}} d_H(\hat \QQ_{n}(\mathbb K(\tau)),\QQ_P(\mathbb  K(\tau)))  \to_{\P^*} 0,  \  \ \sup_{\tau\in \mathcal{T}}d_H(\hat \QQ_{n}( \mathcal K(\tau)),\QQ_P(\mathcal K(\tau)))  & \to_{\P^*} 0.
\end{eqnarray*}

\end{corollary}

The main results  are derived assuming we know the reference distribution $F$ and the associated
depth function $\DD_F$ as well as depth regions $\mathbb{K}(\tau)$ and quantile contours~$\mathcal{K}(\tau)$.  There are cases where these will be approximated numerically
or using data. The same definitions and results extend naturally where these quantities are replaced by uniformly consistent estimators $\hat \DD_{F,n}$, $\hat {\mathbb{K}}_n(\tau)$, and  $\hat {\mathcal{K}}_n(\tau)$:
\begin{equation}\label{depth}
\begin{array}{rc}
\sup_{ u \in K } |  \hat D_{F,n}(u) - D_F(u) | &  \to_{\P^*} 0,   \\
\sup_{\tau\in \mathcal{T}}d_H( \hat {\mathbb K}_n(\tau), \mathbb K(\tau))  & \to_{\P^*} 0, \\
\sup_{\tau\in \mathcal{T}}d_H( \hat {\mathcal K}_n(\tau), \mathcal K(\tau)) &  \to_{\P^*} 0,
\end{array}
\end{equation}
where $K$ is any closed subset of $\mU_0$. These high-level conditions hold trivially for the numerical approximations we use in Section~\ref{sec:comp}. They also hold, for example, for Tukey's halfspace depth under regularity conditions. We will not discuss these conditions here.

\begin{definition}[Empirical MK depth, quantiles, ranks and signs with estimated~$F$]
\label{def:EMK2}Let $F$ be an absolutely continuous reference distribution with support contained in a convex and compact 
region~$\mU \subset \Bbb{R}^d$, and let~$\| \cdot \|$ be a norm on~$\mU$. Given estimators~$\hat \DD_{F,n}$, $\hat {\mathbb{K}}_n(\tau)$ and~$\hat {\mathcal{K}}_n(\tau)$ satisfying (\ref{depth}), the MK empirical quantiles, ranks, signs and depth are defined as follows.
\begin{enumerate}
\item  The {\em MK empirical rank and sign} of $y \in \R^d$ are $\|\hat \RR_n(y)\|$ and~$ \hat \RR_n (y)/\|\hat \RR_P(y)\|$.
\item   The {\em MK empirical} $\tau$-{\em quantile contour} is the set $\hat \QQ_n(\hat{\mK}_n(\tau))$ and the {\em MK empirical depth region} with probability mass $\tau$ is~$\hat \QQ_n(\hat{\bK}_n(\tau))$.
\item The {\em MK empirical depth} of $y \in\R^d$ is the depth of $\hat \RR_n(y)$ under~$\hat {\DD}_{F,n}$:
$$
\hat\DD^{\rm\scriptsize MK}_{P,n}(y):=\hat{\DD}_{F,n}(\hat \RR_n(y)).
$$
\end{enumerate}
\end{definition}

\begin{corollary} \label{cor:eqd2}  Work with conditions of the previous corollary and suppose
that Conditions (\ref{depth}) hold. Then the conclusions of Corollary 3.1 hold and the conclusions
of Corollary 3.2 hold in the following form:
\begin{eqnarray*}\quad \sup_{\tau\in \mathcal{T}} d_H(\hat \QQ_{n}(\hat{\bK}_n(\tau)),\QQ_P(\mathbb  K(\tau)))   \to_{\P^*}   0, \ \ 
\sup_{\tau\in \mathcal{T}}d_H(\hat \QQ_{n}( \hat{\mK}_n(\tau)),\QQ_P(\mathcal K(\tau)))   \to_{\P^*}  0.
\end{eqnarray*}
 \end{corollary}


\section{Computing Empirical Quantiles and Depth Regions}\label{sec:comp}

Here we provide computational characterizations of the empirical quantiles, ranks,
and depth regions for various cases of interest.


\subsection*{Smooth $\hat P_n$ and $\hat F_n$} Suppose $\hat P_n$ and $\hat F_n$ satisfy Caffarelli regularity conditions, so that $\hat\QQ_n=\nabla\hat\psi_n$ and $\hat\RR_n=\nabla\hat\psi^\ast_n$, with $(\hat\psi_n,\hat\psi^\ast_n)$ satisfying (C). The MK empirical vector quantile maps $\hat\QQ_n$ and $\hat\RR_n$ can then be computed with the algorithm of Benamou and Brenier \cite{BB:2000}. 

\subsection*{Discrete $\hat P_n$ and smooth $\hat F_n$}
Suppose now $\hat P_n$ is a discrete estimator of $P$ and~$\hat F_n$~an absolutely continuous distribution with convex compact support $\mU \subset\R^d\!$. 
Let~$\hat P_n$ be of the form $\hat P_n\!=\!\sum_{k=1}^{K_n}p_{k,n}\delta_{y_{k,n}}$ for some integer $K_n$, some nonnegative weights~$p_{1,n},\ldots,p_{K_n,n}$ such that $\sum_{k=1}^{K_n}p_{k,n}=1$,
and $y_{1,n},\ldots,y_{K_n,n}\in\R^d$ . The leading example is when $\hat P_n$ is the empirical distribution of a random sample
$(Y_i)_{i=1}^n$ drawn from $P$.

The MK empirical vector quantile map $\hat\QQ_n$ is then equal (almost everywhere) to the gradient of a convex map $\hat\psi_n$ such that $\nabla\hat\psi_n\#\hat F_n=\hat P_n$, i.e., the $\hat F_n$-almost surely unique map $\hat\QQ_n=\nabla\hat\psi_n$ satisfying the following:
\begin{enumerate}
\item[(1)] $\nabla\hat\psi_n(u)\in\{y_{1,n},\ldots,y_{K_n,n}\}$, for Lebesgue-almost all $u\in\mU$,
\item[(2)] $\hat F_n\left( \{u\in\mU:\;\nabla\hat\psi_n(u)=y_{k,n}  \} \right)=p_{k,n}$, for each $k\in\{1,\ldots,K_n\}$,
\item[(3)] $\hat\psi_n$ is a convex function.
\end{enumerate}
The following characterization of $\hat\psi_n$ specializes Kantorovich duality to this discrete-continuous case (see, e.g.,\cite{EGH}).
\begin{lemma*}
There exist unique (up to an additive constant) weights
$\{v_1^\ast,\ldots,v_n^\ast\}$ such that
$\hat\psi_n(u)=\max_{1\leq k\leq K_n}\{u^\top y_{k,n}-v_k^\ast\}$
satisfies conditions (1),
(2) and (3). The function~$
v\mapsto\int
\hat\psi_nd\hat F_n+\sum_{k=1}^{K_n}p_{k,n}v_k
$
is convex and minimized at
$v^\ast=\{v_1^\ast,\ldots,v_n^\ast\}$.
\end{lemma*}

This lemma allows efficient computation of $\hat\QQ_n$ using a gradient algorithm proposed in \cite{AHA:98}. The map $\hat\psi_n$ is piecewise affine and the empirical vector quantile $\hat\QQ_n$ is piecewise constant. The
correspondence $\hat\QQ_n^{-1}$ defined for each $k\leq K_n$ by
\[
y_{k,n}\mapsto\hat\QQ^{-1}_n(y_{k,n}):=\{u\in\mU:\;\nabla\hat\psi_n(u)=y_{k,n}  \}
\]
maps $\{y_{1,n},\ldots,y_{K_n,n}\}$ into $K_n$ regions of a partition of $\mU$, called a {\em power diagram}.
The estimator $\hat\RR_n$ of the MK vector rank can be computed according to  formula (\ref{eq: empirical Q and R}) after computing the conjugate $\hat \psi^*_n$ of $\hat \psi_n$ via: $\hat \psi^*_n(y) = \sup_{u \in \mU} \{ u^\top y - \hat \psi_n(u) \}.$ The empirical depth,  depth regions, and quantiles can be computed using the depth function, according to their theoretical definitions.

\subsection*{Discrete $\hat P_n$ and $\hat F_n$}
Particularly amenable to computation is the case when both distribution estimators $\hat P_n$ and $\hat F_n$ are discrete with uniformly distributed mass on sets of points of the same cardinality. Let $\hat P_n=\sum_{j=1}^n\delta_{y_j}/n$ for a set $\mathcal Y_n=\{y_1,\ldots,y_n\}$ of points in $\R^d$ and $\hat F_n=\sum_{j=1}^n\delta_{u_j}/n$, for a set $\mathcal U_n=\{u_1,\ldots,u_n\}$ of points in $\R^d$. The restriction of the quantile map $\hat\QQ_n$ to $\mathcal U_n$ is the bijection
$u \longmapsto y=\hat\QQ_n|_{\mathcal U_n}(u)$ from $\mathcal U_n$ onto $\mathcal Y_n$ and $\hat\RR_n|_{\mathcal Y_n}$ is its inverse. The solutions $\hat\QQ_n$ and $\hat\RR_n$ can be computed with any optimal assignment algorithm. More generally, in the case of any two discrete estimators~$\hat P_n$ and~$\hat F_n$, the problem of finding $\hat\QQ_n$ or $\hat\RR_n$ is a linear programming problem.

\subsection*{Visualization of Empirical MK Depth and Quantile Contours}  Whenever $\hat P_n$ is finitely discrete,  then the MK empirical depth regions and quantile contours are finite sets of points. For visualization
purposes it may be helpful to transform them into nicer looking objects which are close to the original objects in terms of Hausdorff distance.   In the example below we used $\alpha$-hulls to create approximations to the depth regions and took the boundaries of the set as a numerical approximation to the quantile contours.    It may also be possible to use polygonization methods such as those in  \cite{DS:88} for $d=2$ and \cite{Gruenbaum:94} for $d=3$.

\begin{example}[Computing MK Depth Regions]   In the example illustrated in Figure~2,
we use a discrete approximation $\hat F_n$ to the spherical uniform reference distribution.  Figure~2 shows the MK empirical depth contours  for the same banana-shaped distribution as in Figure~1. The specific construction to produce Figure~2 is the following:~$\hat P_n$ is the empirical distribution of a random sample $\mathcal Y_n$ drawn from the banana-shaped distribution in $\R^2$, with $n=9999$; $\hat F_n$ is a discrete approximation to $F$ with mass $1/n$ on each of the points in $\mathcal U_n$. The latter is a collection of $99$ evenly spaced points on each of $101$ circles, of evenly spaced radii in $(0,1]$. The sets $\mathcal Y_n$ and $\mathcal U_n$  are matched optimally with the assignment algorithm of the {\em adagio} package in R. MK empirical depth regions are
$\alpha$-hulls of $\hat\QQ_n(\mathcal U_n\cap\mathbb S(\tau))$ for $11$ values of $\tau\in(0,1)$ (see \cite{EKDS:83} for a definition of $\alpha$-hulls). The~$\alpha$-hulls are computed using the {\em alphahull} package in~R, with~$\alpha=0.3$.
The banana-shaped distribution considered is the distribution of the vector~$(X+R\cos\Phi,X^2+R\sin\Phi)$, where~$X$ is uniform on~$[-1,1]$, $\Phi$ is uniform on~$[0,2\pi]$, $Z$ is uniform on $[0,1]$, $X$, $Z$ and~$\Phi$ are independent, and~$R=0.2Z(1+(1-|X|)/2)$. \qed
\end{example}

\appendix

\section{Uniform Convergence of Subdifferentials and Transport Maps}

\subsection{Uniform Convergence of Subdifferentials}
Let $\mathcal{U}$ and $\mathcal{Y}$ be convex, closed subsets of $\R^d$.  A
 pair of convex potentials  $\psi: \mathcal{U} \mapsto \R \cup \{\infty\}$
and $\psi^*: \mathcal{Y} \mapsto \R \cup \{\infty\}$  is a conjugate pair over $(\mU, \mY)$ if, for each $u \in \mU$ and $y \in \mY$,
$$
\psi(u) = \sup_{y \in \mY} [y^\top u - \psi^*(y)], \quad
\psi^*(y) = \sup_{u \in \mU} [y^\top u - \psi(u)].
$$

In the sequel, we consider a fixed  pair $(\psi, \psi^*)$ obeying the following condition.

\begin{itemize}
\item[(C)]
Let $\mathcal{U}$ and $\mathcal{Y}$ be closed, convex subsets of $\R^d$, and $\mU_0 \subset \mU$ and $\mY_0 \subset \mY$ 
some open, non-empty sets in $\R^d$.   Let $\psi: \mathcal{U} \mapsto \R$ and $\psi^*: \mathcal{Y} \mapsto \R$ form a conjugate pair over $(\mU, \mY)$  and possess gradients $\nabla \psi(u)$ for all $u \in \mU_0$, and $\nabla \psi^*(y)$ for all $y \in \mY_0$. The gradients $\nabla \psi |_{\mU_0}:  \mU_0 \mapsto \mY_0$ and $\nabla \psi^* |_{\mY_0}: \mY_0 \mapsto \mU_0$ are homeomorphisms, and~$\nabla \psi|_{\mU_0} = (\nabla \psi^* |_{\mY_0})^{-1}$.
\end{itemize}

We also consider a sequence $(\psi_n, \psi^*_n)$ of conjugate potentials approaching $(\psi, \psi^*)$.

\begin{itemize}
\item[(A)] A sequence of conjugate potentials $(\psi_n, \psi^*_n)$ over $(\mU, \mY)$, with  $n \in \mathbb{N}$,
is such that: $\psi_n(u) \to \psi(u)$ in $\R \cup \{\infty\}$ pointwise in $u$ in a dense subset of $\mU$
and $\psi^*_n(y) \to \psi^*(y)$ in $\R \cup \{\infty\}$  pointwise in $y$ in a dense subset of $\mY$, as $n \to \infty$.
\end{itemize}
Condition (A) is equivalent to requiring that either $\psi_n$ or $\psi^*_n$ converge
pointwise over dense subsets.  There is no loss of generality in stating that both converge.

Define the maps
$$
\QQ(u) := \arg \sup_{y \in \mY}[ y^\top u - \psi^*(y)] , \quad
\RR(y) := \arg \sup_{u \in \mU} [y^\top u - \psi(u)],
$$
for each $u \in \mU_0$ and $y \in \mY_0$. By the envelope theorem,
$$\RR(y) = \nabla \psi^*(y),  \text{ for } y \in \ \mY_0; \ \ \
 \QQ(u) = \nabla \psi(u), \text{ for } u \in  \ \mU_0.$$
Let us define, for each $u \in \mU$ and $y \in \mY$,
\begin{equation}\label{define Qn and Rn}
\QQ_n(u) \in \arg \sup_{y \in \mY} [y^\top u - \psi_n^*(y)] , \quad
\RR_n(y) \in \arg \sup_{u \in \mU} [y^\top u - \psi_n(u)].
\end{equation}

It is useful to note that
$$\RR_n(y) \in \partial \psi^*_n(y) \text{ for }  y \in  \mY;  \ \ \QQ_n(u) \in \partial \psi_n(u) \text{ for } u \in   \mU,$$
where $\partial$ denotes the sub-differential of a convex function; conversely,
any pair of elements of $\partial \psi^*_n(y)$ and $\partial \psi_n(u)$, respectively, could be taken
as solutions to the problem  (\ref{define Qn and Rn}) (by Proposition 2.4 in Villani \cite{villani1}).  Hence, the problem of convergence of
$\QQ_n$ and $\RR_n$ to~$\QQ$ and $\RR$ is equivalent to the problem of convergence of subdifferentials.
Moreover, by Rademacher's theorem,
$\partial \psi^*_n(y) = \{ \nabla \psi_n^*(y) \} $ and~$\partial \psi_n(u) = \{ \nabla \psi_n(u) \}$
almost everywhere with respect to the Lebesgue measure (see, e.g., \cite{villani1}), so the solutions
to (\ref{define Qn and Rn}) are unique almost everywhere on $u \in \mU$ and~$y \in \mY$.

\begin{theorem}[Local uniform convergence of subdifferentials]\label{lemma: ucg}
Suppose that Conditions~(A) and (C) hold.  Then,  as $n \to \infty$, for any compact set $K \subset   
\mU_0$ and any compact set $K' \subset   \mY_0$,
$$
\sup_{ u \in K} \|  \QQ_n(u) - \QQ(u) \|  \to 0,   \quad
\sup_{ y \in K'} \|  \RR_n(y) - \RR(y) \|  \to 0.
$$
\end{theorem}

\begin{remark}
This result appears to be new. It complements the result stated in Lemma 5.4 in Villani \cite{villani:stability} for the case $\mU_0 = \mU=\mY_0 = \mY = \R^d$.  This result also trivially implies convergence in $L^p$ norms, $1 \leq p < \infty$:
$$
\int_{\mU} \| \QQ_n(u) - \QQ(u)\|^p d F(u) \to 0, \quad  \int_{\mY} \| \RR_n(y) - \RR(y)\|^p d P (y)  \to 0,
$$
for probability laws $F$ on $\mU$ and $P$ on $\mY$,  whenever, for some $\bar p>p$, 
$$
\sup_{n \in \mathbb{N}} \int_{\mU}  \| \QQ_n(u)\|^{\bar p} +\|\QQ(u)\|^{p} d F(u) < \infty, \quad
\sup_{n \in \mathbb{N}} \int_{\mY}  \| \RR_n(y)\|^{\bar p} +\|\RR(y)\|^{p} d P(y) < \infty.
$$
Hence, the new result is stronger than available results on convergence in measure (including $L^p$ convergence results) in the optimal transport literature (see, e.g., Villani~\cite{villani1,villani2}).  \qed
\end{remark}

\begin{remark}
 The following example also shows that, in general, our result cannot be strengthened
 to the uniform convergence over entire sets $\mU$ and $\mY$.
Consider the sequence of  potential maps $\psi_n: \mU = [0,1] \mapsto \R$:
$$
\psi_n(u) =  \int_0^u \QQ_n(t) dt, \ \  \QQ_n(t) = t \cdot 1( t \leq 1-1/n) + 10 \cdot 1( t > 1 -1/n).
$$
Then
$
\psi_n(u) = 2^{-1}u^2 1( u \leq 1-1/n) +  \big\{
10(u -(1-1/n)) + 2^{-1}(1-1/n)^2  \big\}1(u> 1- 1/n) $
converges uniformly on $[0,1]$ to
$
\varphi(u) = 2^{-1} u^2.$
The latter potential has the gradient map~$\QQ: [0,1] \mapsto \mY_0 = [0,1]$
defined by $\QQ(t) = t$.  We have that $
\sup_{t \in K}|\QQ_n(t) - \QQ(t)| \to 0
$ for any compact subset $K$ of $(0,1)$. However,  the uniform convergence
over the entire region $[0,1]$ fails, since
$\sup_{t \in [0,1]}|\QQ_n(t) - \QQ(t)| \geq 9$ for all $n$.
Therefore, the theorem cannot be strengthened in general. \qed
\end{remark}

We next consider the behavior of image sets of gradients
defined as follows:
$$
\QQ_{n}(A) := \{\QQ_{n}(u): u \in A\},  \quad \QQ(A) := \{\QQ(u): u \in A\},  \quad A \subseteq K,
$$$$
 \mathrm{R}_{n}(A') := \{ \RR_{n}(y): y \in A'\}, \quad  \RR(A') := \{ \RR(y): y \in A'\}, \quad A' \subseteq K',
$$
where  $K \subset  \  \mU_0$ and  $K' \subset \ \mY_0$ 
are compact sets, and the subsets $A$ and $A'$ are understood to be non-empty.

\begin{corollary}[Convergence of sets of subdifferentials] \label{cor: csg}
Under the conditions of the previous theorem, we have that
$$
\sup_{A  \subseteq K}d_H(\QQ_{n}(A),\QQ(A))  \to 0, \quad
\sup_{A ' \subseteq K'} d_H(\RR_{n}(A'), \RR(A'))  \to 0.
$$
\end{corollary}

\begin{corollary}[Convergence of sets of subdifferentials] \label{cor: csg2}
Assume the conditions of the previous theorem.  For any sequence of sets $\{ A_n\} \subseteq K$ and $\{A'_n \} \subseteq K'$ such that~$d_H(A_n, A) \to 0$ and $d_H(A_n', A')\to 0$ for some sets $A$ and $A'$, we have
$$
d_H(\QQ_{n}(A_n),\QQ(A))  \to 0, \quad
 d_H(\RR_{n}(A'_n), \RR(A'))  \to 0.
$$
\end{corollary}

\subsection{Uniform Convergence of Transport Maps}
We next consider the problem of convergence for potentials
and transport (vector quantile and rank) maps arising from the Kantorovich
dual optimal transport problem.

Equip $\mY$ and $\mU$ with absolutely continuous probability measures $P$  and $F$, respectively,  and let
$$\mY_0 := \text{int}(\text{support}(P)),\qquad \mU_0:= \text{int}(\text{support}(F)).$$
%
We consider  sequences of measures
$P_n$ and $F_n$  approximating $P$ and $F$:
\begin{itemize}
\item[(W)] There are sequences of measures $\{P_n\}_{n \in \mathbb{N}}$ on $\mY$
and $\{F_n\}_{n \in \mathbb{N}}$ on $\mU$, with finite total mass,  that converge
to $P$ and $F$, respectively,  in the topology of weak convergence:
$$d_{\BL} (P_n, P) \to  0, \ \quad d_{\BL} (F_n, F) \to 0 .$$
\end{itemize}

Recall that we defined $\Phi_0(\mU, \mY)$ as a collection of conjugate potentials $(\phi, \phi^*)$
on~$(\mU, \mY)$ such that $\phi(u_0)= 0$ for some fixed point $u_0 \in  \mU_0$.       Let $(\psi_n, \psi^*_n) \in \Phi_0(\mU, \mY)$ solve
the Kantorovich problem for the pair $(P_n,F_n)$:
\begin{equation}\label{empirical kantorovich}
\int \psi_n dF_n + \int \psi_n^* d P_n=  \inf_{ (\phi, \phi^*) \in \Phi_0(\mU, \mY) } \int \varphi dF_n + \int \varphi^* dP_n.
\end{equation}
Also, let $(\psi, \psi^*) \in \Phi_0(\mU, \mY)$ solve the Kantorovich problem  for the pair $(P,F)$:
\begin{equation}\label{kantorovich}
\int \psi dF + \int \psi^* d P = \inf_{ (\phi, \phi^*) \in \Phi_0(\mU, \mY) } \int \varphi dF + \int \varphi^* dP .
\end{equation}
It is known that solutions to these problems  exist; see, e.g., Villani \cite{villani1}.  Recall also that we imposed the normalization
condition in the definition of $\Phi_0(\mU, \mY)$ to pin down the constants.

\begin{theorem}[Local uniform convergence of transport maps]\label{theorem: uct}  Suppose that
 the sets $\mU$ and $\mY$ are compact subsets of $\R^d$, and that the probability measures $P$ and $F$ are absolutely continuous with respect to the Lebesgue measure, with ${\rm support }(P) \subseteq \mY$ and ${\rm support }(F) \subseteq \mU$. Suppose that Condition (W) holds, and that Condition  (C) holds for a solution $(\psi, \psi^*)$ of (\ref{kantorovich}) for the sets $\mU_0$ and  $\mY_0$ defined as above. Then the conclusions of Theorem~\ref{lemma: ucg}
 and Corollary \ref{cor: csg} and Corollary \ref{cor: csg2}. hold.
\end{theorem}


\section{Proofs}

%

\subsection{Proof of Theorem~\ref{lemma: ucg}} The proof  relies on the equivalence of the uniform
and continuous convergence.

\begin{lemma}[Uniform convergence via continuous convergence]
\label{Lemma: Resnick} Let $\mathbb{D}$ and $\mathbb{E}$ be complete
separable metric spaces, with $\mathbb{D}$ compact. Suppose $f: \mathbb{D}
\mapsto\mathbb{E}$ is continuous. Then a sequence of functions $f_{n}:
\mathbb{D} \mapsto\mathbb{E}$ converges to $f$ uniformly on~$\mathbb{D}$ if
and only if, for any convergent sequence $x_{n} \to x$ in $\mathbb{D}$, we
have that $f_{n}(x_{n})\to f(x)$.
\end{lemma}
For the proof, see, e.g.,  Rockafellar and Wets \cite{RW}. The proof also relies on the following convergence result, which is a consequence of Theorem 7.17 in Rockafellar and Wets \cite{RW}.   For a point $a$ and a non-empty set $A$ in $\Bbb{R}^d$,
define $d(a,A) := \inf_{a' \in A} \| a- a'\|$.

\begin{lemma}[Argmin convergence for convex problems]\label{lemma: convex}   Suppose that $g$ is a lower-semi-continuous convex function mapping $\R^d$ to $\R \cup \{+ \infty\}$  that attains a  minimum on the set $\mX_0 =  \arg\inf_{x \in \R^d} g(x)\subset \mathrm{int}(\mD_0)$, where $\mD_0  = \{x \in \R^d: g(x) < \infty \}$, and $ \mathrm{int}(\mD_0)$ is a non-empty, open set in $\R^d$. Let  $\{g_n\}$
be a sequence of convex, lower-semi-continuous functions mapping $\R^d$ to $\R \cup \{+ \infty\}$ and such that
$g_n(x) \to g(x)$ pointwise in $x \in \R_0^d$, where $\R_0^d$ is a countable dense subset of $\R^d$.
Then any $x_n \in \arg\inf_{x \in \R^d} g_n(x)$ obeys
$$
d(x_n, \mX_0) \to 0,
$$
and, in particular, if $\mX_0$ is a singleton $\{x_0\}$, $x_n \to x_0$.
\end{lemma}
The proof of this lemma is given below, immediately after the conclusion of the proof of this theorem.

We define the extension maps $y \mapsto g_{n,u}(y)$ and $u \mapsto \bar g_{n,y}(u)$ mapping $\R^d$ to $\R \cup \{-\infty\}$
$$
g_{n,u} (y) := \left \{ \begin{array}{cc} y^\top u - \psi^*_n(y) & \text{ if  }
y \in \mY \\
- \infty  & \text{ if  }
y \not \in \mY
\end{array} \right.,   \quad
\bar g_{n,y} (u) := \left \{ \begin{array}{cc} y^\top u - \psi_n(u) & \text{ if  }
u \in \mU \\
- \infty  & \text{ if  }
u \not \in \mU.
\end{array} \right.
$$
By the convexity of $\psi_n$ and $\psi^*_n$ over convex, closed sets $\mY$ and $\mU$, we have that  the functions are proper upper-semi-continuous concave functions.  Define the extension maps~$y \mapsto g_{u}(y)$ and $u \mapsto \bar g_{y}(u)$ mapping $\R^d$ to $\R \cup \{-\infty\}$ analogously, by removing the index $n$ above.

Condition (A)  assumes pointwise convergence of $\psi^*_n$ to $\psi^*$
on a dense subset of $\mY$.  By Theorem 7.17 in Rockafellar and Wets \cite{RW}, this implies the uniform convergence of $\psi^*_n$ to
$\psi^*$ on  any compact set $K' \subset \text{ int } \mY$ 
that does not overlap with the boundary of the set $\mD_1 = \{ y \in \mY:  \psi^*(y) < + \infty\}$. Hence, for any sequence~$\{u_n\}$ such that $u_n \to u \in K$, a compact subset of $ \mU_0$, and any $y \in (\text{int }  \mY) \setminus \partial \mD_1$,
  $$g_{n,u_n} (y) = y^\top u_n - \psi^*_n(y) \to g_{u}(y) =  y^\top u - \psi^*(y).$$
Next, consider any $y \not \in \mY$, in which case,  $
g_{n,u_n}(y) = -\infty \to g_u(y) = - \infty.
$ Hence,
$$
g_{n,u_n}(y) \to g_u(y) \text{ in } \R \cup \{ - \infty \} , \text{ for all } y \in \R^d_1 = \R^d \setminus (\partial \mY \cup \partial \mD_1),
$$
where $\R^d_1$ is a dense subset of $\R^d$.  We apply Lemma \ref{lemma: convex} to conclude that
$$
 \arg\sup_{y \in \R^d} g_{n,u_n}(y) \ni \QQ_n(u_n)  \to \QQ(u) \in \arg\sup_{y \in \R^d} g_{u}(y) = \{ \nabla \psi(u) \}.
$$

Take $K$  as any compact subset of $ \mU_0$. The above argument applies for every point~$u \in~\!K$ and every convergent sequence $u_n \to u$. Therefore, since by Assumption~(C) the map $u \mapsto Q(u) = \nabla \psi (u)$ is continuous in $u \in K$, we conclude by the equivalence of the continuous and uniform convergence, Lemma \ref{Lemma: Resnick}, that
$$
\QQ_n(u)  \to \QQ(u) \text{ uniformly in } u \in K.
$$

By symmetry, the proof of the second claim is identical to the proof of the first one.~\!\qed

\subsection{Proof of Lemma \ref{lemma: convex}}  By assumption,
$\mX_0 = \arg \min g \subset \mathrm{int}(\mD_0)$, and $\mX_0$ is convex and closed. Let $x_0$ be an element of $\mX_0$.  We have that, for all $ 0< \eps \leq \eps_0$ with $\eps_0$ such that
$B_{\eps_0}(\mX_0) \subset \mathrm{int}(\mD_0)$,
\begin{equation}\label{g bound}
g(x_0) < \inf_{x \in \partial B_{\eps}(\mX_0)} g(x),
 \end{equation}
where $B_{\eps} (\mX_0) := \{x  \in \R^d: d(x, \mX_0) \leq \eps\}$ is convex and closed.

Fix  an $\eps \in (0, \eps_0]$. By convexity of $g$ and $g_n$ and by Theorem 7.17 in Rockafellar and Wets \cite{RW}, the pointwise convergence of $g_n$ to $g$
on a dense subset of $\R^d$ is equivalent to the uniform convergence of $g_n$ to $g$ on
any compact set $K$ that does not overlap with~$\partial \mD_0$, i.e. $K \cap \partial \mD_0 = \emptyset$.
Hence, $g_n \to g$ uniformly on $B_{\eps_0}(\mX_0)$.   This and (\ref{g bound}) imply that eventually, i.e. for all $n \geq n_\eps$,
$$
g_n(x_0) < \inf_{x \in \partial B_\eps (\mX_0)} g_n(x).
$$
 By convexity of $g_n$, this implies that
$g_n(x_0) < \inf_{x \not \in B_\eps (\mX_0)} g_n(x)$ for all $n \geq n_\eps$,  which is to say that,  for all $n \geq n_\eps$,
$$
\arg\inf g_n = \arg\min g_n  \subset B_\eps (\mX_0).
$$
Since $\eps >0$ can be set as small as desired,  it follows that any $x_n \in \arg \inf g_n$ is such that~$d(x_n, \mX_0) \to 0$. \qed

\subsection{Proof of Corollary \ref{cor: csg}}   By Theorem \ref{lemma: ucg} and the definition of Hausdorff distance, for $A$ denoting non-empty subsets,
\begin{eqnarray*} && \sup_{A \subseteq K}d_H(\QQ_{n}(A), \QQ(A)) \\
&&  = \sup_{A \subseteq K}  \left ( \sup_{u \in A}\inf_{\bar u \in A}\| \QQ_n(\bar u) - \QQ(u)  \|
  \vee  \sup_{\bar u\in A} \inf_{u \in A} \| \QQ_n(\bar u) - \QQ(u)  \| \right) \\
&&  \leq \sup_{A \subseteq K} \left(
 \sup_{u \in A} \| \QQ_n(u) - \QQ(u)  \|  \vee  \sup_{\bar u \in A} \| \QQ_n(\bar u) - \QQ(\bar u)  \|
 \right) \\
&& =   \sup_{ u \in K} \|  \QQ_n(u) - \QQ(u) \|  \to 0.
 \end{eqnarray*}
The proof of the second claim is identical. \qed

\subsection{Proof of Corollary \ref{cor: csg2}}   We have that
\begin{eqnarray*}
 d_H(\QQ_n(A_n), \QQ(A)) \!\!\!\! \!&\!\leq\!&\!\! \! d_H(\QQ_n(A_n), \QQ(A_n)) + d_H(\QQ(A_n), \QQ(A))\\ \! \!&\leq\!&\! \!\! \sup_{A \subseteq K} \!d_H(\QQ_n(A), Q(A)) +\!  \sup_{\bar u, u \in K}\! \{\QQ(\bar u) - \QQ(u)\!:  \|\bar u - u\| \leq d_H(A_n, A) \!\} \\ &\to& 0,
\end{eqnarray*}
where the first inequality holds by the triangle inequality, the second inequality holds
by definition and by $A_n, A \subseteq K$, and the last conclusion follows by  Corollary \ref{cor: csg} and continuity
of the map $u \longmapsto Q(u)$ on $u \in K$.

The proof of the second claim is identical. \qed

\subsection{Proof of Theorem~\ref{theorem: uct}}

\textsc{Step 1.}   Here we show that the set of conjugate pairs is compact in the topology of uniform convergence.
First we notice that, for any pair~$(\phi, \phi^*) \in \Phi_0(\mU, \mY)$,
 \begin{eqnarray*}
 \| \varphi \|_{\BL(\mU)} & \leq &  (2 \|\mY\| \|\mU\|) \vee \|\mY\| < \infty, \ \
 \| \varphi^* \|_{\BL(\mY)}  \leq   (2 \|\mY\| \|\mU\|) \vee \|\mU\| < \infty,
 \end{eqnarray*}
 with $\|A\| := \sup_{a \in A} \|a\|$ for $A \subseteq \R^d$, where we have used the fact that $\phi(u_0) = 0$ for some $u_0 \in \mU$ as well as  compactness of $\mY$ and $\mU$.

The Arzela-Ascoli Theorem implies that $\Phi_0(\mU, \mY)$  is relatively compact in the topology of the uniform convergence. We want to show compactness, namely that this set is also closed. For this we need to show that
all uniformly convergent subsequences $(\phi_n, \phi_n^*)_{n \in \mathbb{N}'}$  (where $\mathbb{N}' \subseteq \mathbb{N}$) have the limit point in this set: 
$$
(\phi, \phi^*) := \lim_{n \in \mathbb{N}'} (\phi_n, \phi_n^*)  \in \Phi_0( \mU,  \mY).
$$
This is true, since uniform limits of convex functions are necessarily convex (see \cite{RW}), and since
\begin{eqnarray*}
\phi(u) & = & \lim_{n \in \mathbb{N}'} \left [   \sup_{y \in \mY} [u^\top y - \phi^*_n(y)] \right ] \\
 & \leq & \limsup_{n \in \mathbb{N}'} \left [  \sup_{y \in \mY} [u^\top y - \phi^*(y)] + \sup_{y \in \mY} | \phi^*_n(y) - \phi^*(y) |  \right ]  =  \sup_{y \in \mY} [u^\top y - \varphi^*(y)];
\end{eqnarray*}
and \begin{eqnarray*}
\phi(u) & = & \lim_{n \in \mathbb{N}'} \left [   \sup_{y \in \mY} [u^\top y - \phi^*_n(y)] \right ] \\
 & \geq & \liminf_{n \in \mathbb{N}'} \left [  \sup_{y \in \mY} [u^\top y - \phi^*(y)] - \sup_{y \in \mY} | \phi^*_n(y) - \phi^*(y) |  \right ]  =  \sup_{y \in \mY} [u^\top y - \varphi^*(y)].
\end{eqnarray*}
Analogously, $\varphi^*(y) = \sup_{u \in \mU} [u^\top y - \varphi(y)]$.\vspace{2mm}

\textsc{ Step 2.}  The claim here is that
\begin{equation}\label{value}
I_n: = \int \psi_n d F_n + \int \psi_n^* d P_n \to_{n \in \mathbb{N}} \int \psi dF + \int \psi^* dP =: I_0.
 \end{equation}
Indeed,
$$
I_n \leq \int \psi d F_n + \int \psi^* d P_n   \to_{n \in \mathbb{N}} I_0,
$$
where the inequality holds by definition, and the convergence holds by
$$ \left | \int \psi d (F_n - F) \right |  +     \left |\int \psi^* d(P_n - P) \right  |  \lesssim d_{\BL} (F_n, F) + d_{\BL} (P_n, P)  \to 0,$$
where $x \lesssim y$ means $x \leq A y$,  for some constant $A$ that does not depend on $n$. Moreover, by definition, 
$$
II_n:= \int \psi_n dF + \int \psi^*_n d P \geq I_0,
$$
 but  $$ \left | I_n - II_n\right  | \leq  \left | \int \psi_n d (F_n - F)\right  | + \left | \int \psi_n^* d(P_n - P) \right |
 \lesssim d_{\BL} (F_n, F) + d_{\BL} (P_n, P)  \to 0.$$

\textsc{Step 3.} Here we conclude.

First, we observe that the solution pair $(\psi, \psi^*)$ to the limit Kantorovich problem is unique on $\mU_0\times\mY_0$ in the sense that any other
solution $(\phi, \phi^*)$ agrees with $(\psi, \psi^*)$ on~$\mU_0\times\mY_0$.  Indeed,
suppose that $\phi(u_1) \neq \psi(u_1)$ for some $ u_1 \in \mU_0$.
By  the uniform continuity of elements of $\Phi_0(\mU, \mY)$ and openness of $\mU_0$, there exists a ball $B_\eps(u_1) \subset \mU_0$ such that $\psi(u) \neq \phi(u)$ for all $ u \in B_\eps(u_1)$. By  the normalization assumption $\phi(u_0) =  \psi(u_0)=0$,
there does not exist a constant $c \neq 0$ such that
$\psi(u) = \phi(u) + c$ for all~$u \in \mU_0$, so this must imply that $\nabla \psi (u) \neq \nabla \phi(u)$ on a set $K \subset \mU_0$ of positive measure (otherwise, if they disagree only on a set of measure zero, we would have
$\psi(u) - \psi(u_0) = \int_0^1 \nabla \psi(u_0 + v^{\top}(u-u_0))^\top (u-u_0) d v = \int_0^1 \nabla \phi(u_0 + v^{\top}(u-u_0)) ^\top(u-u_0) d v
= \phi(u)- \phi(u_0)$ for almost all~$u \in B_\eps(u_1)$, which is a contradiction). However,  the statement $\nabla \psi \neq \nabla \phi$ on a set~$K \subset \mU_0$ of positive Lebesgue measure would
contradict the fact that any solution $\psi$ or $\phi$ of the Kantorovich problem must obey $$
\int h \circ \nabla \phi  d F =  \int h \circ \nabla \psi d F = \int  h dP ,
$$
for each bounded continuous $h$, i.e. that $\nabla \phi \# F = \nabla \psi \#  F = P$, established on p.72 in Villani \cite{villani1}.
Analogous arguments apply to establish uniqueness of $\psi^*$ on the set $\mY_0$.

Second, we can split $\mathbb{N}$ into subsequences $\mathbb{N} = \cup_{j=1}^\infty \mathbb{N}_{j}$
such that,  for each $j$,
\begin{equation}\label{unif}
(\psi_n, \psi^*_n) \to_{n \in \mathbb{N}_j} (\varphi_j, \varphi^*_j) \in \Phi_0( \mU, \mY), \ \text{ uniformly on $\mU \times \mY$. }
\end{equation}
But by Step 2 this means that
$$
\int \varphi_j d F + \int \varphi^*_j dP = \int \psi dF + \int \psi^* dP.
$$
It must be that  each pair $(\phi_j, \phi_j^*)$ is the solution to the limit Kantorovich problem,
and by the uniqueness established above we have that
$$
(\phi_j, \phi_j^*) = (\psi, \psi^*) \text{ on } \mU_0 \times \mY_0.
$$
By Condition (C) we have that, for $u \in \mU_0$ and $y \in \mY_0$:
$$
\QQ(u) = \nabla \psi(u) = \nabla \phi_j(u), \quad  \RR(u) = \nabla \psi^*(u) = \nabla \varphi^*_j(u).
$$
By (\ref{unif}) and Condition (C) we can invoke Theorem~\ref{lemma: ucg} to conclude that
$\QQ_n \to \QQ$ uniformly on compact subsets of $\mU_0$ and $\RR_n \to \RR$ uniformly
on compact subsets of~$\mY_0$. \qed

\subsection{Proof of Lemma \ref{ergodic}} The proof is a variant of standard arguments,
for example, those given in (\cite{PRW:1999}, proof of Theorem~2.1), so is delegated to the Supplemental Appendix.  \qed

\subsection{Proof of Theorem~\ref{theorem: empirical}}  The proof is an immediate consequence of the Extended Continuous Mapping Theorem, as given in van der Vaart and Wellner \cite{VvW}, Theorem~\ref{lemma: ucg} and Corollary \ref{cor: csg}.

The theorem, specialized to our context, reads as follows:  Let $\mathbb{D}$ and $\mathbb{E} $ be normed spaces
and let $x \in \mathbb{D}$. Let $\mathbb{D}_n \subseteq \mathbb{D}$ be arbitrary subsets and $g_n: \mathbb{D}_n \mapsto \mathbb{E}$ be arbitrary maps~($n \geq 0$), such that for every sequence $x_n \in
\mathbb{D}_n$ such that $x_{n} \to x$, along a subsequence, we have that $g_{n}(x_{n}) \to g_0(x)$, along the same subsequence. Then, for arbitrary
(i.e. possibly non-measurable) maps $X_n: \Omega \mapsto \mathbb{D}_n$  such that  $X_n \to_{\P^*} x$, we have that~$g_n(X_n) \to_{\P^*} g_0(x)$.

In our case $X_n = (\hat P_n, \hat F_n)$ is a stochastic element of $\mathbb{D}$,  viewed as an arbitrary map from $\Omega$ to $\mathbb{D}$,
and $x = (P, F)$ is a non-stochastic element of $\mathbb{D}$, where $\mathbb{D}$
is  the space of linear operators $\mathbb{D}$ acting on the space of bounded Lipschitz functions.
This space can be equipped with the norm (see notation section):
\[
\|\cdot\|_{\mathbb{D}}:
 \| (x_1, x_2) \|_{\mathbb{D}} =  \| x_1 \|_{\BL(\mY)} \vee \| x_2 \|_{\BL(\mU)}.
\]  Moreover,
$X_n \to_{\Pr^*} x$ with respect to this norm, i.e.
$$
\| X_n - x\|_{\mathbb{D}}:= \| \hat P_n - P \|_{\BL(\mY)} \vee \|\hat F_n-F\|_{\BL(\mU)} \to_{\Pr^*} 0.
$$
 Then $g_n(X_n) := (\hat \QQ_n, \hat \RR_n)$ and $g (x) := (\QQ, \RR)$ are
viewed as elements of the space\linebreak  $\mathbb{E}= \ell^\infty (K\times K', \R^d \times \R^d)$ 
of bounded functions  mapping $K\times K'$ to $\R^d \times \R^d$, equipped with the supremum
norm. The maps have the continuity property: if $\|x_n - x\|_{\mathbb{D}} \to 0$
along a subsequence, then $\|g_n(x_n) - g(x)\|_{\mathbb{E}} \to 0$ along the same subsequence,
as established by Theorem~\ref{lemma: ucg}.  Hence
conclude that  $g_n(X_n) \to_{\P^*} g(x)$.

The second claim follows by the Extended Continuous
Mapping Theorem and  Corollary \ref{cor: csg}. \qed

\subsection{Proof of Corollaries \ref{cor:ERS}, \ref{cor:eqd}, and \ref{cor:eqd2}}
 Corollaries \ref{cor:ERS} and \ref{cor:eqd} follow  by Theorem~\ref{theorem: empirical}  and the Extended Continuous Mapping Theorem;
and Corollary \ref{cor:eqd2} follows by Theorem~\ref{theorem: empirical}, the Extended Continuous Mapping Theorem, and Corollary ~\ref{cor: csg2}.   \qed

\bibliographystyle{plain}

\end{document}